\documentclass{amsart}
\usepackage{amssymb}
\usepackage{eucal}
\theoremstyle{plain}
\newtheorem{theorem}{Theorem}
\newtheorem{corollary}{Corollary}

\newtheorem{lemma}{Lemma}

\theoremstyle{definition}

\theoremstyle{remark}

\numberwithin{equation}{section}

\newcommand{\R}{\mathbb R}

\newcommand{\Z}{\mathbb Z}
\newcommand{\C}{\mathbb C}

\newcommand{\Rn}{\mathbb R^n}

\begin{document}

\title[Unique continuation]{Uniqueness Properties of Solutions to 
Schr\"odinger Equations}
\author{L. Escauriaza}
\address[L. Escauriaza]{UPV/EHU\\Dpto. de Matem\'aticas\\Apto. 644, 48080 Bilbao,
Spain.}
\email{luis.escauriaza@ehu.es}
\thanks{The first and fourth authors are supported  by MEC grant,
MTM2004-03029, the second and third authors by NSF grants DMS-0968472 and
DMS-0800967 respectively}
\author{C. E. Kenig}
\address[C. E. Kenig]{Department of Mathematics\\University of Chicago\\Chicago, Il.
60637 \\USA.}
\email{cek@math.uchicago.edu}
\author{G. Ponce}
\address[G. Ponce]{Department of Mathematics\\ University of California\\ Santa
Barbara, CA 93106\\ USA.}
\email{ponce@math.ucsb.edu}
\author{L. Vega}
\address[L. Vega]{UPV/EHU\\Dpto. de Matem\'aticas\\Apto. 644, 48080 Bilbao, Spain.}
\email{luis.vega@ehu.es}

\maketitle


\section{Introduction}\label{S: Introduction}
To place the subject of this paper in perspective, we start out with a brief discussion 
of  unique continuation. Consider solutions to

\begin{equation}
\label{aa1}
\Delta u(x)=\sum_{j=1}^n\frac{\partial^2u}{\partial x_j^2}(x)=0,
\end{equation}
(harmonic functions) in the unit ball $\,\{x\in \R^n\,:\,|x|<1\}$. When $n=2$, these functions are real parts
of holomorphic functions, and so, if they vanish of infinite order at $x=0$, they must vanish identically. We call this the strong unique continuation property (s.u.c.p.).  The same result holds for $n>2$, since harmonic functions are still real analytic in  $\,\{x\in \R^n\,:\,|x|<1\}$. In fact,  it is well-known that if $\,P(x,D)$ is a linear elliptic differential operator with real analytic coefficients, and $\,P(x,D)u=0$ in a open set $\,\Omega\subset \R^n$, then $u$  is real analytic in $\,\Omega$. Hence, the (s.u.c.p.) also holds for such solutions.  Through the work of Hadamard \cite{Had} on the uniqueness of the Cauchy problem (which is closely related to the strong unique continuation property discussed earlier) it became clear (for applications in nonlinear problems) that it would be desirable to establish the strong unique continuation property for operators whose coefficients are not necessarily real analytic, or even $\, C^{\infty}$. The first results in this direction were found in the pioneering work of Carleman \cite{Car} (when $n=2$) and M\"uller \cite{Mu} (when $n>2$), who proved the (s.u.c.p) for 
$$
P(x,D)=\Delta+V(x),\;\;\;\;\;\text{with}\;\;\;\;\;\;V\in L^{\infty}_{loc}(\R^n). 
$$
 In order to establish his result, Carleman introduced a method (the method of \lq\lq Carleman estimates") which has permeated the subject ever since. In this context, an example of a Carleman estimate is :  
 
\vskip.05in

\emph{ For $\,f\in C^{\infty}_0(\{x\in \R^n\,:\,|x|<1\}-\{0\})$, $\,\alpha>0$ 
 and  
 $$
 w(r)=r\,\exp (\,\int_0^r\frac{e^{-s}-1}{s} ds),
 $$
 one has 
  \begin{equation}
 \label{1aa}
 \alpha^3\,\int\,w^{-1-2\alpha}(|x|) f^2(x) dx\leq c\,\int w^{2-2\alpha}(|x|)\,|\Delta f(x)|^2 dx,
 \end{equation}
 with $\,c\,$ independent of $\,\alpha$} 
 \vskip.05in

 For a proof of this estimate, see \cite{EsVe}, \cite{BoKe}. The (s.u.c.p.) of Carleman-M\"uller follows easily from \eqref{1aa} 
 (see \cite{Ke2} for instance).  
 
 In the late 1950's and 1960's there was a great deal of activity on the subject of (s.u.c.p.) and the closely related uniqueness 
 in the Cauchy problem, some highlights being \cite{AKS} and \cite{Cal} respectively, both of which use the method of Carleman estimates. These results and methods 
 have had a multitude of applications to many areas of analysis, including to non-linear problems. (For a recent example, see \cite{KeMe} for an application to energy critical
 non-linear wave equations).
 
  In connection with the Carleman-M\"uller (s.u.c.p.) a natural question is :  How fast is a solution $u$ allowed to vanish, before it must vanish identically?
  
  By considering $n=2$, $u(x_1,x_2)=\Re (x_1+ix_2)^N$, we see that to make sense of the question, a normalization is required, for instance
  $$
  \sup _{|x|<3/4} |u(x)|\geq 1,\;\;\;\;\;\;\;\|u\|_{L^{\infty}(|x|<1)}<\infty.
  $$
We refer to questions of this type as \lq\lq quantitative unique continuation". It is also of interest to consider unique continuation type questions around the point at infinity. For instance, a conjecture of E. M. Landis \cite{KoLa} was :  if
$$
\Delta u +V u =0,\;\;\;\;x\in\R^n,\;\;\;\text{with}\;\;\;\|V\|_{\infty}\leq 1 ,\;\;\;\|u\|_{\infty}<\infty,
$$
and for some $\,\epsilon>0$  one has 
$$
|u(x)|\leq c_{\epsilon}\,e^{-c_{\epsilon}|x|^{1+\epsilon}},
$$
then  $\,u\equiv 0$.  

For the case of complex valued potentials $V(x)$,
 this conjecture was disproved by Meshkov \cite {Me} who constructed $V,\,u,\,u\not \equiv 0$ with 
 $$
 |u(x)|\leq c\, e^{-c|x|^{4/3}},\;\;\;\;n\geq 2.
 $$
 Meshkov also showed that if
 $$
 |u(x)|\leq c_{\epsilon}\, e^{-c_{\epsilon}|x|^{4/3+\epsilon}},\;\;\;\;\text{for some}\;\;\;\;\epsilon>0,
 $$
 then $\,u\equiv 0$. 
 
 It turns out that a \lq\lq quantitative" formulation of this can also be proved, as it was done in \cite{BoKe}, and this was crucial for the resolution in \cite{BoKe} of a long-standing problem in disordered media, namely Anderson localization near the bottom of the spectrum, for the continuous Anderson-Bernoulli model in $\,\R^n,\,n\geq 1$.
 
 Next, we turn to versions of unique continuation for evolution equations. We start with parabolic equations and consider solutions of
 $$
 \partial_t u-\Delta u=W\cdot \nabla u+ Vu,\;\;\;\;\;\;\;\text{with}\;\;\;\;\;\;\|W\|_{\infty}+\|V\|_{\infty} <\infty,
 $$
 (or equivalently $|\partial_tu-\Delta u|\leq M(|\nabla u|+|u|)$).
 Using a parabolic analog of the Carleman estimate described earlier, one can show that if
 $$
|\partial_tu-\Delta u|\leq M(|\nabla u|+|u|),\;\;\;\;\;\;(x,t)\in \{x\in\R^n:|x|<4R\}\times [t_0,t_1],\;\;\;R>0,
$$
with $\,|u(x)|\leq A$ and
$$
u\equiv 0,\;\;\;(x,t)\in \{x\in\R^n:R<|x|<4R\}\times [t_0,t_1],
$$
then
$$
u\equiv 0,\;\;\;(x,t)\in \{x\in\R^n:|x|<R\}\times [t_0,t_1].
$$

We call this type of result \lq\lq unique continuation through spatial boundaries", (see \cite{EsVe}, \cite{Ve} and references therein for this type of result 
and strengthenings of it). This result is closely related to the \lq\lq elliptic" (s.u.c.p.) discussed before. On the other hand, for parabolic equations, there is also a 
\lq\lq backward uniqueness" principle, which is very useful in applications to control theory (see \cite{lm60} for an early result in this direction) : Consider solutions to
$$
 |\partial_tu-\Delta u|\leq M(|\nabla u|+|u|),\;\;\;\;\;(x,t)\in\R^n\times (0,1],
 $$
with $\,\|u\|_{\infty}\leq A$. Then, if $\,u(\cdot,1)\equiv 0$, we must have $\,u\equiv 0$. 
This result is also proved through Carleman estimates (see \cite{lm60}). 

Recently, a strengthening of this result has been obtained in \cite{EsSS}, where one considers solutions only defined in $\,R^n_{+}\times (0,1]$, $\,R^n_{+}=\{(x_1,..,x_n)\in\R^n\,:\,x_1>0\}$, without any assumptions on $\,u\,$ at $\,x_1=0$, and still obtains the \lq\lq backward uniqueness" result. This   strengthening had an important application to non-linear equations, allowing the authors of \cite{EsSS}  to establish a long-standing conjecture of J. Leray on regularity and uniqueness of solutions to the Navier-Stokes equations (see also \cite{Se} for a recent extension).

Finally, we turn to dispersive equations. Typical examples of these are the $k$-generalized KdV equation
\begin{equation}
\label{kdv1}
 \partial_t u +\partial_x^3 u+ u^k\partial_xu=0,\;\;\;\;\;\; (x,t)\in \R\times
\R,\;\,\;k\in \Z^+,
\end{equation}
and the non-linear Schr\"odinger equation
\begin{equation}
\label{ae1}
 \partial_t u =i( \Delta u \pm |u|^{p-1}u),\;\;\;\;\;\; (x,t)\in
\mathbb{R}^n\times \R,\;\;\,p>1.
\end{equation}
These equations model phenomena of wave propagation and have been extensively studied in the last 30 years or so. 

For these equations,\lq\lq unique continuation through spatial boundaries " also holds, as it was shown by Saut-Scheurer \cite{SaSc} for the KdV-type equations and by Izakov \cite{Iza} for Shr\"odinger type equations. (All of these results were established trough Carleman estimates). These equations  however are time reversible (no preferred time direction) and so \lq\lq backward uniqueness" is immediate, unlike in parabolic problems.
Once more in connection with control theory, this time for dispersive equations, Zhang \cite{BZ} showed, for solutions of 
\begin{equation} 
\label{zhang} 
 \partial_t u =i( \partial_x^2 u \pm |u|^{2}u),\;\;\;\;\;\; (x,t)\in
\mathbb{R}\times[0,1],
\end{equation} 
that if $\,u(x,t)=0$ for $(x,t)\in (-\infty,a)\times\{0, 1\}$ (or $(x,t)\in (a,\infty)\times\{0, 1\}$) for some $a\in\R$, the $\,u\equiv 0$.
Zhang's proof was based on the inverse scattering method which uses that this is a completely integrable model, and did not apply to
other non-linearities or dimensions. This type of result  was extended to the $k$-generalized KdV \eqref{kdv1} and the general non-linear Schr\"odinger equation in \eqref{ae1} in all dimensions 
(where inverse scattering is no longer available) using suitable Carleman estimates (see \cite{KPV02}, \cite{IK04}, \cite{IK06}, and references therein). 

For recent surveys of the results presented so far, see \cite{Ke1}, \cite{Ke2}.

 Returning to \lq\lq backward uniqueness" for parabolic equations, in analogy with Landis' \lq\lq elliptic" conjecture mentioned earlier, Landis-Oleinik \cite{LaOl}  conjectured that in the  \lq\lq backward uniqueness"  result one can replace the hypothesis $\,u(\cdot,1)\equiv 0$ with the weaker one
 $$
 |u(x,1)|\leq c_{\epsilon}\,e^{- c_{\epsilon}|x|^{2+\epsilon}},\;\;\;\text{for some }\;\;\;\epsilon>0.
 $$
This is indeed true and was established in \cite{EKPV06a} and \cite{Ng}. Similarly, one can conjecture (as it was done in \cite{EKPV08b}) that for Schr\"odinger equations, if 
 $$
|u(x,0)|+ |u(x,1)|\leq c_{\epsilon}\,e^{- c_{\epsilon}|x|^{2+\epsilon}},\;\;\;\text{for some }\;\;\;\epsilon>0,
 $$
then $\,u\equiv 0$. This was established in \cite{EKPV06a}.

In analogy with the improvement of   \lq\lq backward uniqueness"  in \cite{EsSS}, one can show that it suffices to deal with solutions in $\,\R^n_{+}\times(0,1]$ 
(for parabolic problems) and require
$$
|u(x,1)|\leq c_{\epsilon}\,e^{- c_{\epsilon}x_1^{2+\epsilon}},\;\;\;x_1>0,\;\;\;\;\text{for some }\;\;\;\epsilon>0,
$$
to conclude that $\,u\equiv 0$ (\cite{Ng}), and that for the Schr\"odinger equations it suffices to have $\,u\,$ a solution in $\,\R^n_{+}\times[0,1]$, with
$$
|u(x,0)|+ |u(x,1)|\leq c_{\epsilon}\,e^{- c_{\epsilon}x_1^{2+\epsilon}},\;\;\;x_1>0,\;\;\;\;\text{for some }\;\;\;\epsilon>0,
$$
to conclude that $\,u\equiv 0$, as we will prove in section 5 of this paper. 

In \cite{EKPV06} it was pointed out for the first time (see also \cite{Cha})  that both the results in \cite{EKPV06a} and in \cite{EKPV06}, in the case of the free heat 
equation
$$
\partial_tu=\Delta u,
$$
and the free Schr\"odinger equation
$$
\partial_tu=i \Delta u,
$$
respectively,  are in fact a corollary of the more precise Hardy uncertainty principle for the Fourier transform, which says :

\vskip.05in

\emph{If $f(x)=O(e^{-|x|^2/\beta^2})$, $\widehat f(\xi)=O(e^{-4|\xi|^2/\alpha^2})$ and 
$1/\alpha\beta>1/4$, then $f\equiv 0$, and if $1/\alpha\beta=1/4$, $f(x)=ce^{-|x|^2/\beta^2}$}  as will be discussed below.
\vskip.05in

Thus, in a series of papers (\cite{EKPV06}-\cite{EKPV10}, \cite{CEKPV}) we took up the task of finding the sharp version of the Hardy uncertainty principle, in the 
context of evolution equations. The results obtained have already yielded new results on non-linear equations. For instance in \cite{EKPV08m} and \cite{EKPV10}
we have found applications to the decay of concentration profiles of possible self-similar type blow-up  solutions of non-linear Schr\"odnger equations and to the decay of possible solitary wave type solutions of non-linear Schr\"odinger equations.

In the rest of this work we shall review some of our recent results concerning unique continuation
properties of
solutions of Schr\"odinger equations of the form
\begin{equation}
\label{e1}
 \partial_t u =i( \Delta u + F(x,t,u,\bar u)),\;\;\;\;\;\; (x,t)\in
\mathbb{R}^n\times \R.
\end{equation}

We shall be  mainly interested in the case where 
\begin{equation}
\label{F1a}
F(x,t,u,\bar u)=V(x,t) u(x,t) 
\end{equation}
is describing the evolution of the Schr\"odinger flow with a time dependent 
potential $V(x,t)$, and 
in the semi-linear case
\begin{equation}
\label{F1b}
F(x,t,u,\bar u)= F(u,\bar u), 
\end{equation}
with $ F: \C \times \C\to \C$, $F(0,0)=\partial_uF(0,0)=\partial_{\bar u}F(0,0)=0$. 

Let us consider a familiar dispersive model, 
the $k$-generalized
Korteweg-de Vries equation \eqref{kdv1} and recall a theorem established in \cite{EKPV07} :
 
\begin{theorem}
\label{theorem1}
 There exists $c_0>0$ such that for any pair  
 $$
u_1,\,u_2\in C([0,1]:H^4(R)\cap L^2(|x|^2dx))
$$
of  solutions of   
\eqref{kdv1} 
 such that if 
\begin{equation}
 \label{3:2}
 u_1(\cdot,0)-u_2(\cdot,0),\,\;\, u_1(\cdot,1)-u_2(\cdot,1)\in
L^2(e^{c_0x_{+}^{3/2}}dx),
\end{equation}  then $u_1\equiv u_2$.
  \end{theorem}
 
 Above we have used the notation:   $ x_{+}=max\{x;\,0\}$. 
 
 Notice that taking $u_2\equiv 0\,$ Theorem \ref{theorem1} gives a restriction on
the possible decay 
of a non-trivial solution of 
\eqref{kdv1}
at two different times.
The power $3/2$ in the exponent in \eqref{3:2} reflects the asymptotic behavior of
the Airy function. 
More precisely, the solution 
of the initial value problem (IVP) 
\begin{equation}
 \begin{aligned}
 \begin{cases}
 \partial_t v + \partial_x^3 v=0,\\
 v(x,0)=v_0(x),
 \end{cases}
 \end{aligned}
\end{equation}
is given by the group $\{U(t)\,:\,t\in R\}$
 $$
 U(t)v_0(x)=\frac{1}{\root{3}\of{3t}}\,Ai\left(\frac{\cdot}{\root
{3}\of{3t}}\right)\ast v_0(x), 
 $$
where
 $$
Ai(x)=c\,\int_{-\infty}^{\infty}\,e^{ ix\xi+i \xi^3 }\,d\xi,
 $$
 is the Airy function which satisfies the estimate 
$$
 |Ai(x)|\leq c (1+x_{-})^{-1/4}\,e^{-c x_{+}^{3/2}}.
 $$
 
It was also shown in \cite{EKPV07} that Theorem 1 is optimal :
\begin{theorem}
\label{theorem2}
There exists $ \,u_0\in S(\R),\;u_0\not \equiv 0$  and $\Delta T>0$ such that 
 the IVP associated 
 to the k-gKdV equation \eqref{kdv1} 
 with data $u_0$ has solution 
  $$
  u\in C([0,\Delta T] :  \mathbb S(\R)),
  $$
  satisfying 
$$
 |u(x,t)|\leq \tilde d \,e^{-x^{3/2}/3}, \;\;\;\;\;\;\;\;x>1,\;\;\;\,t\in [0,\Delta T],
 $$
for some constant $\tilde d>0$. 
 \end{theorem}

In the case of the free Schr\"odinger group $\{e^{it\Delta}\,:\,t\in\R\}$ 
$$
e^{it\Delta}u_0(x)=(e^{-i|\xi|^2t} \widehat
u_0)^\lor(x)=\frac{e^{i|\cdot|^2/4t}}{(4\pi i t)^{n/2}}*u_0(x),
$$
the fundamental solution does not decay.
However, one has the identity 
\begin{equation}
 \label{formula1}
\begin{aligned}
&u(x,t)= e^{it\Delta}u_0(x)= \int_{\R^n} \frac{e^{i|x-y|^2/4t}}{(4\pi i t)^{n/2}}\,
u_0(y)\,dy\\
\\
&=\frac{e^{i|x|^2/4t}}{(4\pi i t)^{n/2}} \int_{\R^n}e^{-2ix\cdot y/4t} e^{i|y|^2/4t}
u_0(y)\,dy\\
\\
&= \frac{e^{i|x|^2/4t}}{(2 i t)^{n/2}}\;
\widehat{\;(e^{i|\cdot|^2/4t}u_0)\,}\left(\frac{x}{2 t}\right),
\end{aligned}
\end{equation}
where
$$
\widehat f(\xi)=(2\pi)^{-n/2} \int_{\R^n} e^{-i\xi\cdot x} f(x)dx.
$$

Hence, 
$$
c_t e^{-i|x|^2/4t} \,u(x,t) = \widehat{(e^{i|\cdot|^2/4t}u_0)}\left(\frac{x}{2
t}\right),\,\,\,\,\,\,\,c_t=(2 i t)^{n/2},
$$
which tells us that 
$e^{-i|x|^2/4t} \,u(x,t)$  is a multiple of the rescaled Fourier transform
 of $\;e^{i|y|^2/4t}u_0(y)$. Thus, as we pointed out earlier, the behavior of the solution of the free
Schr\"odinger equation is closely 
related to uncertainty principles for the Fourier transform.
We shall study  
these uncertainty principles
and their relation with the uniqueness properties of the solution of the
Schr\"odinger equation \eqref{e1}. 
In the early $1930$'s N. Wiener's remark (see \cite{Hardy}, \cite{In}, and \cite{Mo}):
\vskip.1in
``a pair of transforms $f$ and $g$ ($\widehat f$) cannot 
both be very small'', 
\vskip.1in
\noindent motivated the works of G. H. Hardy \cite{Hardy}, G. W. Morgan \cite{Mo}, and A. E. Ingham \cite{In} which
will 
be considered in detail in this note.  
However, before that we shall return to a review of some previous results concerning uniqueness
properties of solutions
of the Schr\"odinger equation which we mentioned earlier and which were not motivated by the formula 
\eqref{formula1}.

 For solutions $u(x,t)$ of the $1$-D cubic Schr\"odinger equation \eqref{zhang} 
B. Y. Zhang \cite{BZ} showed :

\vskip.05in
\emph{If  $u(x,t)=0$ for $(x,t)\in (-\infty, a)\times \{0,1\}\,\, ($or $(x,t)\in
(a,\infty)\times \{0,1\})\,$ for some 
$\,a\in \R$, then $u\equiv 0$.}
\vskip.03in
 As it was mentioned before, his 
proof is based on the inverse scattering method, which uses the fact that the
equation in \eqref{zhang} is a completely integrable model. 
\vskip.05in
In \cite{KPV02} it was proved under general assumptions on $F$
 in \eqref{F1b} that :
\vskip.03in
\emph{If $u_1,\,u_2\in C([0,1]:H^s(\R^n))$, with $\,s>\max \{n/2; \,2\}\, $ are
solutions of the equation \eqref{e1} with
$F$ as in \eqref{F1b} such that
$$
u_1(x,t)=u_2(x,t),\;\;\;(x,t)\in \Gamma^c_{x_0}\times \{0,1\},
$$
where $ \Gamma^c_{x_0}$ denotes the complement of a cone  $\Gamma_{x_0}$ with vertex
$x_0\in \R^n$ and opening $<180^0$, 
then
$u_1\equiv u_2$.}

(For further results in this direction see \cite{KPV02}, \cite{IK04}, \cite{IK06},
and 
references therein). 
\vskip.03in
A key step in the proof in \cite{KPV02} was the following uniform exponential decay
estimate:

 \begin{lemma}\label{ultimo} There exists $\epsilon_0>0$ such that if
\begin{equation}
\label{hyp2}
\mathbb V:\mathbb R^n\times [0,1]\to\mathbb C,\;\;\;\;\text{with}\;\;\;\;
\|\mathbb V\|_{L^1_tL^{\infty}_x}\leq \epsilon_0,
\end{equation}
and $u\in C([0,1]:L^2(\mathbb R^n))$ is a strong solution of the IVP
\begin{equation}
\begin{cases}
\begin{aligned}
\label{eq1}
&\partial_tu=i(\Delta +\mathbb V(x,t))u+\mathbb G(x,t),\\
&u(x,0)=u_0(x),
\end{aligned}
\end{cases}
\end{equation}
with
\begin{equation}
\label{hyp3} u_0,\,u_1\equiv u(\,\cdot\,,1)\in
L^2(e^{2\lambda\cdot x}dx),\;\mathbb G\in L^1([0,1]:L^2(e^{2\lambda\cdot
x}dx)),
\end{equation}
for some $\lambda\in\mathbb R^n$, then there exists $c_n$ independent of
$\lambda$ such that
\begin{equation}
\begin{aligned}
\label{uno}
&\sup_{0\leq t\leq 1}\| e^{\lambda\cdot x} u(\,\cdot\,,t)\|_{L^2(\mathbb \R^n)} \\
&\leq c_n
\Big(\|e^{\lambda\cdot x} u_0\|_{L^2(\mathbb \R^n)} + \|e^{\lambda\cdot x}
u_1\|_{L^2(\mathbb \R^n)} +\int_0^1
\|e^{\lambda\cdot x}\, \mathbb G(\cdot, t)\|_{L^2(\mathbb \R^n)} dt\Big).
\end{aligned}
\end{equation}
\end{lemma}
\vspace{0,07 cm}
Notice that in the above result one assumes the existence of a reference
$L^2$-solution $u$ of the equation \eqref{eq1} and then under the 
hypotheses \eqref{hyp2} and \eqref{hyp3} shows that the exponential decay in the
time interval $[0,1]$ is preserved.

The estimate \eqref{uno} can be combined with the subordination formula 
\begin{equation}
\label{est1}
e^{\gamma |x|^p/p}\simeq \int_{\R^n} \,e^{\gamma^{1/p}\lambda\cdot
x-|\lambda|^q/q}\,
|\lambda|^{n(q-2)/2}\,d\lambda,\,\,\,\forall\, x\in \R^n\,\,\,\text{and}\,\,\,p>1,
\end{equation}
to get that for any $\alpha>0$ and $ a>1$
\vspace{0,07 cm}
\begin{equation}
\label{dos}
\begin{aligned}
\sup_{0\leq t\leq 1}\| e^{\alpha|x|^a} u(\,\cdot\,,t)\|_{L^2(\mathbb \R^n)} &\\
\leq c_n
\Big(\|e^{\alpha|x|^a} u_0\|_{L^2(\mathbb \R^n)} +& \|e^{\alpha|x|^a}
u_1\|_{L^2(\mathbb \R^n)} +\int_0^1
\|e^{\alpha|x|^a}\, \mathbb G(\cdot, t)\|_{L^2(\mathbb \R^n)} dt\Big).
\end{aligned}
\end{equation}
 
 Under appropriate assumptions on the potential $V(x,t)$ in \eqref{F1a} one writes
$$
V(x,t)u=  \chi_{R} V(x,t)u + (1-\chi_{R}) V(x,t)u = \mathbb V(x,t)u + \mathbb G(x,t),
$$
with $\chi_R\in C^{\infty}_0,\,$ $\chi_R(x)=1,\,|x|<R$, supported in $|x|<2R$,  and
applies the estimate 
\eqref{dos} by fixing $\,R\,$ sufficiently large. Also under appropriate hypothesis
on $F$ and $u$  a similar argument 
can be used for 
the semi-linear equation in 
\eqref{F1b}.
\vskip.03in
 The estimate \eqref{dos} gives a control on the decay of the solution in the whole
time interval 
in terms of that  at the end  points and that of  the \lq\lq external force''. As we shall see
below a key  idea  
will be to get improvements of this estimate based on  logarithmically convex
versions of it.

We recall that if one considers the equation \eqref{e1} with initial data $u_0\in
\mathbb S(\R^n)$
and a smooth potential $V(x,t)$ in \eqref{F1a} 
or smooth nonlinearity $F$ in \eqref{F1b}, it follows that the corresponding
solution satisfies that
$u\in C([-T,T] :\mathbb S(\R^n))$. This can be proved using the commutative property
of the operators
$$
L=\partial_t-i\Delta,\;\;\;\;\;\;\text{and}\;\;\;\;\;\;\Gamma_j=x_j+2t\partial_{x_j},\,\,\,j=1,..,n,
$$
see \cite{HKT1}-\cite{HKT2}. From the proof of this fact one also has that the
persistence property of the solution $u=u(x,t)$ 
(i.e. if the data $u_0\in X$, a function space, then the corresponding solution
$u(\cdot)$ describes 
a continuous curve in $X$, $u\in C([-T,T]:X)$, $\,T>0$)
with data $u_0\in L^2(|x|^m)$
can only hold if $u_0\in H^s(\R^n)$ with $s\geq 2m$.  
Roughly speaking, for exponential weights one has a more involved argument where 
 the time direction plays a role. 
Considering the IVP for the one dimensional free 
Schr\"odinger equation 
\begin{equation}
\label{*}
\begin{cases}
\begin{aligned}
&\partial_t u=i\partial_x^2u,\;\;\;\;\;\;\;\;\;\;\;\,x,\,t \in \R,\\
&u(x,0)=u_0(x)\in L^2(\R),
\end{aligned}
\end{cases}
\end{equation}
and assuming that $ e^{\beta x}u_0 \in L^2(\R),\,\,\beta>0$, then one formally has
that \
$$
v(x,t)=e^{\beta x}u(x,t)
$$
satisfies the equation
$$
\partial_t v=i(\partial_x-\beta)^2v.
$$
Thus, 
$$
v(x,\pm 1)=e^{\beta x}u(x,\pm 1)\in L^2(\R)\,\;\;\;\text{if}\;\;\;\,e^{\pm 2\beta \xi}\,
\widehat{e^{\beta x}u_0}\in L^2(\R).
$$
However, if we knew that 
$ e^{\beta x}u(x,1),\;\;e^{\beta x}u(x,-1) \in L^2(\R)$ integrating forward in time 
the positive frequencies of $e^{\beta x}u(x,t)$ and backward in time the negative
frequencies of $e^{\beta x}u(x,t)$
one gets an estimate similar to that in \eqref{uno} with $\lambda=\beta$ and
$\mathbb G=0$. This argument  motivates the
idea behind Lemma \ref{ultimo} and its proof. 

\vskip.2in

 The rest of this paper is organized as follows: section 2 contains the results
related to Hardy's uncertainty principle including a short discussion on the version of this  principle 
in terms of   the heat flow. 
Section 3 those concerned with Morgan's uncertainty principle. In section 4 we
shall consider the limiting case  in section 3. Also, section 4 includes the statements 
of some related forthcoming results. 
Earlier in the introduction we have discussed uniqueness results obtained under the assumption that the
solution vanishes 
at two different time in a semi-space
(see \cite{BZ}, \cite{IK04}, \cite{IK06}, \cite{EKPV08b}). In section 2 similar uniqueness results will be
established under a Gaussian decay hypothesis, 
in the whole  space. In section 5 we shall obtain a unifying  result, i.e. a
uniqueness result under  Gaussian decay in 
a semi-space of $\,\R^n$ at two different times. The appendix contains an abstract
lemma and a corollary which will be used in
the previous sections.

\section{Hardy's Uncertainty Principle}\label{hardy}
In \cite{Hardy} G. H. Hardy's  proved the following one dimensional ($n=1$) result: 
\vskip.07in
\emph{If 
$f(x)=O(e^{-|x|^2/\beta^2})$, $\widehat f(\xi)=O(e^{-4|\xi|^2/\alpha^2})$ and 
$1/\alpha\beta>1/4$, then $f\equiv 0$. 
\newline
Also, if $1/\alpha\beta=1/4$, $f(x)$ is a constant multiple of $e^{-|x|^2/\beta^2}$.}

\vskip.06in
To our knowledge  the available proofs of this result and its variants use complex
analysis, mainly appropriate versions 
of the  
Phragm\'en-Lindel\"of principle. 
There has also been considerable interest in a better understanding of this result
and on extensions of it 
to other settings: 
 \cite{bonamie1}, \cite{bonamie2}, \cite{CoPr}, \cite{Ho}, and \cite{SST}. In
particular, the extension of 
Hardy's result to higher 
dimension $n\geq 2$ (via Radon transform) was given in \cite{SST}.

The formula \eqref{formula1} allows us to re-write this uncertainty principle in terms of the
solution of the IVP 
 for the free  Schr\"odinger equation 
$$
\begin{cases}
\begin{aligned}
&\partial_tu=i\triangle u, \;\;\,\,(x,t)\in\R^n\times (0,+\infty),\\
&u(x,0)=u_0(x),
\end{aligned}
\end{cases}
$$
in the following manner : 
\vskip.05in

\emph{If $u(x,0)=O(e^{-|x|^2/\beta^2})$, $u(x,T)=O(e^{-|x|^2/\alpha^2})$ and
$T/\alpha\beta> 1/4$, then 
$u\equiv 0$. Also, if $T/\alpha\beta=1/4$, $u$ has as initial data $u_0$ equal to a
constant multiple of 
$e^{-\left(1/\beta^2+i/4T\right)|y|^2}$.}
\vskip.05in

The corresponding $L^2$-version of Hardy's uncertainty principle was established in
\cite{CoPr2} : 

\vskip.05in
\emph{If $\,e^{|x|^2/\beta^2}f$, $\,e^{4|\xi |^2/\alpha^2}\widehat f$ are in
$L^2(\R^n)$ and 
$1/\alpha\beta\ge 1/ 4$, then $f\equiv 0$.}

\vskip.05in
In terms of the solution of the Schr\"odinger equation it states :

\vskip.05in
\emph{If $\,e^{|x|^2/\beta^2}u(x,0)$, $\,e^{|\xi |^2/\alpha^2}u(x,T)$ are in
$L^2(\R^n)$ and 
$T/\alpha\beta\ge 1/4$, then $u\equiv 0$.}
\vskip.05in

More generally, it was shown in \cite{CoPr2} that :

\vskip.05in
\emph{If $\,e^{|x|^2/\beta^2}f\in L^p(\R^n)$, 
$\,e^{4|\xi |^2/\alpha^2}\widehat f\in L^q(\R^n)$, 
$p, q\in [1,\infty]\,$ with at least 
one  of them finite  and $1/\alpha\beta\ge 1/ 4$, then $f\equiv 0$.}
\vskip.05in

In \cite{EKPV08b} we proved a uniqueness result for solutions of \eqref{e1} with $F$ as
in \eqref{F1a} for bounded potentials 
$V$ verifying that either, 
$$
V(x,t)=V_1(x)+V_2(x,t),
$$
with $V_1$ real-valued and
$$
\sup_{[0,T]}\|e^{T^2|x|^2/\left(\alpha
t+\beta\left(T-t\right)\right)^2}V_2(t)\|_{L^\infty(\R^n)}<+\infty,
$$
or
\begin{equation}
 \label{condition}
\lim_{R\rightarrow +\infty}\int_0^T\|V(t)\|_{L^\infty(\R^n\setminus B_R)}\,dt =0.
\end{equation}

More precisely, it was shown  that the only solution $u\in C([0,T], L^2(\R^n))$ to
\eqref{e1} with $F=V(x,t)u$, 
verifying 
\begin{equation}
\label{E: condicion fundamental}
\|e^{|x|^2/\beta^2}u(0)\|_{L^2(\R^n)}+\|e^{|x|^2/\alpha^2}u(T)\|_{L^2(\R^n)}<+\infty
\end{equation}
with $\,T/\alpha\beta>1/ 2$ and $V$ satisfying one of the above conditions
is the zero solution.
Notice that this result differs by a factor of $1/2$ from that for the solution of
the free
 Schr\"odinger equation given by the $L^2$-version of the Hardy uncertainty principle
described above ($T/\alpha\beta\ge 1/4$).

In \cite{EKPV09} we showed that the  
optimal version of Hardy's uncertainty principle in terms of $L^2$-norms, as
established in \cite{CoPr2}, 
holds for solutions of
\begin{equation}
\label{E: 1.11}
\partial_tu=i\left(\triangle u+V(x,t)u\right), \,\,\,\, (x,t)\in \R^n\times [0,T],
\end{equation}
such that \eqref{E: condicion fundamental} holds with $T/\alpha\beta>1/4$ and for
many general bounded potentials $V(x,t)$, 
while it fails for some complex-valued potentials in the end-point case,
$T/\alpha\beta=1/4$. 

\vskip.05in
\begin{theorem}\label{T: hardytimeindepent}
Let $u\in C([0,T]):L^2(\R^n))$ be a solution of the equation
\eqref{E: 1.11}. If there exist positive constants
$\alpha$ and $\beta$ such that  $T/\alpha\beta > 1/4$, 
and
$$
\|e^{|x|^2/\beta^2}u(0)\|_{L^2(\R^n)},\,\,\,\,\|e^{|x|^2/\alpha^2}u(T)\|_{L^2(\R^n)}<\infty,
$$
and the potential $V$ is  bounded and either, 
$V(x,t)=V_1(x)+V_2(x,t)$, with $V_1$ real-valued and 
$$
\sup_{[0,T]}\|e^{T^2|x|^2/\left(\alpha t+\beta \left(T-t\right)\right)^2}V_2(t)
\|_{L^\infty(\R^n)} < +\infty
$$
or 
$$
\lim_{R\rightarrow +\infty}\|V\|_{L^1([0,T], L^\infty(\R^n\setminus B_R)}=0.
$$
 Then, $u\equiv 0$.
\end{theorem}

We remark that there are no assumptions on the size of the potential in the given
class
or on the dimension and that we do not assume any decay of the gradient, neither of
the solutions or 
of the time-independent potential or any \it{a priori  }\rm regularity on this
potential or the solution. 

\vskip.03in

\begin{theorem}\label{T: hardytimeindepent2}
Assume that $T/\alpha\beta=1/4$. Then, there is a smooth complex-valued potential
$V$ verifying 
$$
|V(x,t)|\lesssim\frac 1{1+|x|^2},\, (x,t)\in  \R^n\times [0,T],
$$
 and a nonzero smooth function $u\in C^\infty([0,T],\mathcal S(\R^n))$ solution of
\eqref{E: 1.11}
such that  
\begin{equation}
 \label{007}
\|e^{|x|^2/\beta^2}u(0)\|_{L^2(\R^n)},\,\,\,\,\|e^{|x|^2/\alpha^2}u(T)\|_{L^2(\R^n)}<\infty.
\end{equation} 
\end{theorem}

\vskip.03in

Our proof of Theorem \ref{T: hardytimeindepent} does not use any complex analysis,
giving,  in particular, a new proof (up to the end-point) of the $L^2$-version of 
Hardy's uncertainty principle 
for the Fourier transform. It is based on Carleman estimates for certain evolutions.
More precisely, it is based on the 
convexity and log-convexity properties present for the solutions of these evolutions.
Thus, the 
convexity and log-convexity of appropriate $L^2$-quantities play the role of the
Phragm\'en-Lindel\"of principle.
We observe that the product of log-convex functions is log-convex  which, roughly
speaking, replaces the fact that the 
product of analytic functions is analytic. 

In \cite{CEKPV} in collaboration with M. Cowling, we gave new proofs, based
only on \it{real variable  }\rm techniques, of both the $L^2$-version of the Hardy
uncertainty
principle and the original Hardy's uncertainty
principle $ (L^{\infty}$) $n$-dimensional version for the Fourier transform as stated at the beginning of this
section, 
including the end point  case $1/\alpha \,\beta=1/4$.

 Returning to Theorem \ref{T: hardytimeindepent}  as 
 a by product of our  proof, we obtain the following optimal 
interior estimate for 
the Gaussian decay of solutions to \eqref{E: 1.11}.

\begin{theorem}\label{T: lamejora}
Assume that $\,u\,$ and $\,V\,$ verify the hypothesis in Theorem \ref{T:
hardytimeindepent} 
and $\,T/\alpha\beta\le 1/4$. Then,

\begin{equation}
 \label{oda}
\begin{aligned}
&\sup_{[0,T]}\|e^{a(t)|x|^2}u(t)\|_{L^2(\R^n)} +
\| \sqrt{t(T-t)}\nabla \left(e^{\left(a(t)+\frac{i\dot
a(t)}{8a(t)}\right)|x|^2}u\right)\|_{L^2(\R^n\times [0,T])}\\
&\le N\left[\|e^{|x|^2/\beta^2}u(0)\|_{L^2(\R^n)}+
\|e^{|x|^2/\alpha^2}u(T)\|_{L^2(\R^n)}\right],
\end{aligned}
\end{equation}
where
\[a(t)=\frac {\alpha\beta RT}{2\left(\alpha t+\beta (T-t)\right)^2+2R^2\left(\alpha
t - \beta (T-t)\right)^2}\ ,\]
$R$ is the smallest root of the equation 
$$
\frac T{\alpha\beta}=\frac R{2\left(1+R^2\right)}
$$
and $N$ depends on $T$, $\alpha$, $\beta$ and the conditions on the potential $V$ in
Theorem \ref{T: hardytimeindepent}.
\end{theorem} 

One has  that $1/a(t)$ is convex and attains its minimum value in the interior of
$[0,T]$, when
$$
|\alpha-\beta|<R^2\left(\alpha+\beta\right).
$$

To see the optimality of Theorem \ref{T: lamejora}, we write
\begin{equation}
\label{E: el enemigo}
u_R(x,t)=R^{-\frac n2}\left(t-\frac iR\right)^{-\frac n2}e^{-\frac{|x|^2}{4i(t-\frac
iR)}}=
\left(Rt-i\right)^{-\frac n2}e^{-\frac{(R-iR^2t)}{4(1+R^2t^2)}\,|x|^2},
\end{equation}
which is a free wave (i.e. $V\equiv 0$, in \eqref{E: 1.11}) satisfying in $\R^n\times
[-1,1]$ the corresponding time 
translated conditions in Theorem \ref{T: lamejora} with $T=2$ and
$$
\frac1{\beta^2}=\frac1{\alpha^2}=\mu=\frac R{4\left(1+R^2\right)}\le\frac 18\, .
$$
Moreover
$$
\frac R{4\left(1+R^2t^2\right)}\, ,
$$
is increasing in the $R$-variable, when  $0<R\le 1$ and $-1\le t\le 1$.

Our improvement over the results in \cite{EKPV06} and \cite{EKPV08b} is a consequence 
of the possibility of extending 
the following argument (for the case of free waves) to prove Theorem \ref{T:
hardytimeindepent} 
(a non-free wave case).
 
We recall the conformal or Appell transformation: If   $u(y,s)$ verifies
\begin{equation}
\label{2.1}
\partial_su=i\left(\triangle
u+V(y,s)u+F(y,s)\right),\;\;\;\;\;\;\;(y,s)\in \R^n\times [0,1],
\end{equation}
and $\alpha$ and $\beta$ are positive,  then
\begin{equation}
\label{2.2}
\widetilde u(x,t)=\left(\tfrac{\sqrt{\alpha\beta}}{\alpha(1-t)+\beta
t}\right)^{\frac n2}u\left(\tfrac{\sqrt{\alpha\beta}\,
x}{\alpha(1-t)+\beta t}, \tfrac{\beta t}{\alpha(1-t)+\beta
t}\right)e^{\frac{\left(\alpha-\beta\right) |x|^2}{4i(\alpha(1-t)+\beta
t)}},
\end{equation}
verifies
\begin{equation}
\label{2.3}
\partial_t\widetilde u=i\left(\triangle \widetilde u+\widetilde
V(x,t)\widetilde u+\widetilde F(x,t)\right),\;\; \text{in}\  \R^n\times
[0,1],
\end{equation}
with
\begin{equation}
\label{potencial}
\widetilde V(x,t)=\tfrac{\alpha\beta}{\left(\alpha(1-t)+\beta
t\right)^2}\,V\left(\tfrac{\sqrt{\alpha\beta}\, x}{\alpha(1-t)+\beta t},
\tfrac{\beta t}{\alpha(1-t)+\beta t}\right),
\end{equation}
and
\begin{equation}
\label{externalforce}
\widetilde F(x,t)=\left(\tfrac{\sqrt{\alpha\beta}}{\alpha(1-t)+\beta
t}\right)^{\frac n2+2}F\left(\tfrac{\sqrt{\alpha\beta}\,
x}{\alpha(1-t)+\beta t}, \tfrac{\beta t}{\alpha(1-t)+\beta
t}\right)e^{\frac{\left(\alpha-\beta\right) |x|^2}{4i(\alpha(1-t)+\beta
t)}}.
\end{equation}
Thus,
 to prove Theorem \ref{T: hardytimeindepent} for free waves, it suffices to consider $u\in
C([-1,1], L^2(R^n))$ 
being a solution of 
\begin{equation}
\label{E: free wave}
\partial_tu-=i\triangle u,\,\,\,\,(x,t)\in R\times [-1,1],
\end{equation}
and
\begin{equation}
\label{E: decaimineto}
\|e^{\mu |x|^2}u(-1)\|_{L^2(R^n)}+\|e^{\mu |x|^2}u(1)\|_{L^2(\R^n)}<+\infty,
\end{equation}
for some $\mu >0$. 

 The main idea consists of showing that either $u\equiv 0$ or there is a function 
$\theta_{R}: [-1,1]\longrightarrow [0,1]$ such that
\begin{equation}
\label{E: gaussian improvement}
\|e^{\frac{R|x|^2}{4\left(1+R^2t^2\right)}}u(t)\|_{L^2(R^n)}\le 
\|e^{\mu |x|^2}u(-1)\|_{L^2(R^n)}^{\theta_{R}(t)}\|e^{\mu
|x|^2}u(1)\|_{L^2(\R^n)}^{1-\theta_{R}(t)},
\end{equation}
where $R$ is the smallest root of the equation
$$
\mu =\frac{R}{4\left(1+R^2\right)}\ .
$$
This gives the optimal improvement of the Gaussian decay of a free wave verifying
\eqref{E: decaimineto}  
and we also see that if $\mu > 1/8$, then $u$ is  zero. 

The proof of these facts relies on new logarithmic convexity properties of free
waves verifying 
\eqref{E: decaimineto} and 
on those already established in \cite{EKPV08b}. In \cite[Theorem 3]{EKPV08b}, the
positivity of the 
space-time commutator 
of the symmetric and skew-symmetric parts of the operator,
$$
e^{\mu |x|^2}\left(\partial_t-i\triangle\right)e^{-\mu |x|^2},
$$
is used to prove that  $\|e^{\mu |x|^2}u(t)\|_{L^2(\R^n)}$ is logarithmically convex
in $[-1,1]$. More precisely,
defining 
$$
f(x,t) = e^{\mu |x|^2}u(x,t)=e^{it\Delta}u_0(x),
$$
it follows that
$$
e^{\mu |x|^2}\left(\partial_t-i\triangle\right)u=e^{\mu
|x|^2}\left(\partial_t-i\triangle\right)(e^{-\mu |x|^2}f)
=\partial_t f -\mathcal S f-\mathcal A f,
$$
where $\mathcal S$ is symmetric and $\mathcal A$ skew-symmetric with
$$
\mathcal S=- i\mu(4 \,x\cdot \nabla + 2n),\,\;\,\;\;\;\;\mathcal A=i(\Delta+4\mu^2
\,|x|^2),
$$
so that
$$
[\mathcal S;\mathcal A] = - 8 \mu (\nabla\cdot I \nabla) + 16 \mu^2\,|x|^2.
$$

 Formally, using the abstract Lemma \ref{L: freq1} (see  the appendix) and the
Heisenberg inequality 
$$
\|f\|^2_{L^2(\R^n)}\leq \frac{2}{n} \,\|\,|x|f\|_{L^2(\R^n)}\,\|\,\nabla f\|_{L^2(\R^n)},
$$ 
whose proof follows by integration by parts, one sees that 
$$
H(t)=\|f(t)\|^2_{L^2(\R^n)}=\|e^{\mu |x|^2}u(t)\|_{L^2(\R^n)}
$$
is logarithmically convex so 
$$
\|e^{\mu |x|^2}u(t)\|_{L^2(\R^n)}\le \|e^{\mu |x|^2} u(-1)\|_{L^2(R^n)}^{\frac{1-t}2}
\|e^{\mu |x|^2} u(1)\|_{L^2(\R^n)}^{\frac{1+t}2},
$$
when, $-1\le t\le 1$.

 Setting $a_1\equiv \mu$,  we begin an iterative  process, 
where at the $k$-th step, we have $k$ smooth even functions, 
$a_j:[-1,1]\longrightarrow (0,+\infty)$, $1\le j\le k$, such that
$$
\mu\equiv a_1<a_2<\dots<a_k\in  (-1,1),
$$
$$
F(a_i)> 0,\ a_j(1)=\mu,\ j=1,\dots,k,
$$
where
$$
F(a)=\frac 1a\left(\ddot a-\frac{3\dot a^2}{2a\,}+32a^3\right)
$$
and functions $\theta_j:[-1,1]\longrightarrow [0,1]$, $1\le j\le k$, such that for
$t\in [.1,1]$
\begin{equation}
\label{E: algoagradable}
\|e^{a_j(t) |x|^2}u(t)\|_{L^2(R^n)}\le 
\|e^{\mu |x|^2}u(-1)\|_{L^2(R^n)}^{\theta_j(t)}\|e^{\mu
|x|^2}u(1)\|_{L^2(\R^n)}^{1-\theta_j(t)}.
\end{equation}
 
These estimates follow from the construction of the functions $a_i$, while the
method strongly 
relies on the following formal 
convexity properties of free waves: 
\begin{equation}
\label{E: algo fundamental}
\partial_t\left(\frac 1a\partial_t\log{H_b}\right)\ge -\frac{2\ddot b^2|\xi|^2}{F(a)},
\end{equation}
\begin{equation}
\label{E: el control del gradiente}
\partial_t\left(\frac 1a\partial_tH\right)\ge
\epsilon_a\int_{R^n}e^{a|x|^2}\left(|\nabla u|^2+|x|^2|u|^2\right)\,dx,
\end{equation}
where 
$$
H_b(t)=\|e^{a(t)|x+ b(t)\xi|^2}u(t)\|_{L^2(R^n)}^2\ ,\ 
H(t)=\|e^{a(t)|x|^2}u(t)\|_{L^2(\R^n)}^2,
$$
$\xi\in \R^n$ and $a, b: [-1,1]\longrightarrow R$ are smooth functions with 
$$
a> 0,\quad \;\;\;\;\;F(a)>0 \;\;\;\;\;\text{in}\;\;\;\;[-1,1].
$$

Once the $k$-th step is completed, we take  $a=a_k$ in \eqref{E: algo fundamental}
with a certain choice of 
$b=b_k$, verifying $b(-1)=b(1)=0$ and then, a certain test is performed. When the
answer to the test is positive, 
it follows that $u\equiv 0$. Otherwise, the logarithmic convexity associated to
\eqref{E: algo fundamental} allows 
us to find a new smooth function $a_{k+1}$ in $[-1,1]$ with
$$
a_1<a_2<\dots<a_k<a_{k+1},\,\, (-1,1),
$$
and verifying the same properties as $\,a_1,\dots,a_k$.

When the process is infinite, we have \eqref{E: algoagradable} for all $k\ge 1$ and
there are two possibilities:
$$
\text{ either }\,\,\,\,\,\,\,
\lim_{k\to +\infty}a_k(0)=+\infty,\,\,\,\,\,\,\, \text{or }\,\,\,\,\,\,\,\lim_{k\to
+\infty}a_k(0)<+\infty.
$$
 In the first case and \eqref{E: algoagradable} 
one has  that $u\equiv 0$, while in the second, the sequence $a_k$ is shown to
converge to an even function
 $a$ verifying
\begin{equation}
\label{ole}
\begin{cases}
\ddot a-\frac{3\dot a^2}{2a\,\,}+32a^3=0,\,\,\,\,\, [-1,1]         \\
a(1)=\mu.
\end{cases}
\end{equation}
Because
$$
a(t)=\frac{R}{4\left(1+R^2t^2\right)},\,\,\, \,\quad R\in \R^+,
$$
are all the possible even solutions of this equation, $a$ must be one of them and 
$$
\mu =\frac{R}{4\left(1+R^2\right)},
$$
for some $R>0$. In particular, $u\equiv 0$, when $\mu >1/8$.

\vskip.05in

As it was already mentioned above, our  proof of Theorem \ref{T: hardytimeindepent}
(the case of
non-zero potentials $V=V(x,t)$),  
is based on the extension of the above convexity properties to the non-free case.
\vskip.05in

 Theorem \ref{T: hardytimeindepent2} establishes the sharpness of the result in
Theorem \ref{T: hardytimeindepent} 
by giving an example of a  complex valued potential $V(x,t)$ 
verifying \eqref{condition} and a  non-trivial solution $u\in C([0,T]:L^2(\R^n))$ 
of \eqref{E: 1.11} for which  
\eqref{E: condicion fundamental} holds with $T/\alpha\beta =1/4$.
Thus, one may ask  :  Is it possible to construct 
a real valued potential $V(x,t)$ verifying the same properties, i.e. 
satisfying \eqref{condition}
and having a non-trivial solution $u\in C([0,T]:L^2(\R^n))$ of \eqref{E: 1.11} such
that 
\eqref{E: condicion fundamental} holds with $T/\alpha\beta =1/4\,$?  

The same question concerning the sharpness of the above result presents itself in
the case 
of time independent potentials $V=V(x)$. In this regard, we consider the stationary
problem 
\begin{equation}
\label{estatic}
 \Delta w + V(x) w =0,\,\,\,\, x\in \R^n,\,\,V\in  L^{\infty}(\R^n),
\end{equation}
and recall  V. Z. Meshkov's result in \cite{Me} : 
\vskip.07in
\emph{ If $w\in H^2_{loc}(\R^n)$ is a solution of \eqref{estatic} such that
\begin{equation}
 \label{43}
\int_{\R^n} e^{a|x|^{4/3}}|w(x)|^2dx<\infty,\,\,\,\,\forall a>0,
\end{equation}
then $\,u\equiv 0$.}
\vskip.04in
Moreover, it was also proved in \cite{Me}  that for complex potentials $V$, the
exponent $4/3$ in 
\eqref{43} is optimal. However, it has been conjectured that for real valued potentials
the optimal exponent should be 1, (see also \cite{BoKe} for a quantitative form of these results and applications to Anderson localization of
Bernoulli models). 
\vskip.05in

More generally,  it was established in \cite{EKPV10}, (see also \cite{cruz}) : 
\vskip.07in
\emph{ If $w\in H^2_{loc}(\R^n)$ is a solution of \eqref{estatic} 
with
a complex valued potential $V$ satisfying
$$
V(x)=V_1(x)+V_2(x),
$$
such that 
\begin{equation}
 \label{123}
|V_1(x)|\leq \frac{c_1}{(1+|x|^2)^{\alpha/2}},\,\,\,\,\alpha\in [0,1/2),
\end{equation}
and $V_2$ real valued supported in $\,\{x\,:\,|x|\geq 1\}$ 
such that
$$
-(\partial_r V_2(x))^- < \frac{c_2}{|x|^{2\alpha}},\,\,\,\,a^-=\min\{a;0\}.
$$
Then there exists $a=a(\|V\|_{\infty};c_1;c_2;\alpha)>0$ such that if
\begin{equation}
 \label{43a}
\int_{\R^n} e^{a|x|^{r}}|w(x)|^2dx<\infty,\,\,\,\,\,r=(4-2\alpha)/3,
\end{equation}
then $\,u\equiv 0$.}
\vskip.05in
In addition, one can take the value $r=1$ in \eqref{43} by assuming  $\alpha>1/2$ in
\eqref{123}.
\vskip.05in
It was also proved in \cite{cruz}  that for complex potentials these results for
$\alpha\in [0,1/2)$ are sharp. 
\vskip.05in
By noticing  that given a solution $\phi(x) $ of the eigenvalue problem
\begin{equation}
\label{eigen}
 \Delta \phi + \widetilde V(x) \phi =\lambda \phi,\,\,\,\, x\in \R^n,
\end{equation}
with $\lambda \in\R, $ then $V(x)=\widetilde V(x)+\lambda$ satisfies the hypothesis 
of the previous result
and
$$
u(x,t)= e^{it\lambda} \,\phi(x),
$$
solves the evolution equation
\begin{equation}
\label{evo-notime}
 \partial_t u=i(\Delta u + V(x)u),\,\,\,\, x\in \R^n,\,t\in\R,
\end{equation}
one gets a lower bound for the value of the strongest possible decay rate of  non-trivial
solutions 
$u(x,t)$ of \eqref{evo-notime} at two different times.

\vskip.05in

As a direct  consequence of  Theorem \ref{T: hardytimeindepent} we have the
following application concerning 
the uniqueness of solutions for semi-linear equations of the form \eqref{e1} with
$F$ as in \eqref{F1b}. 

\begin{theorem}
\label{Theorem NL2}

Let $u_1$ and $u_2$ be strong solutions in $C([0,T],H^k(\R^n)), \,k>n/2$ of the
equation \eqref{e1} with
$F$ as in \eqref{F1b}  such that $\,F\in C^k$  and $F(0)=\partial_uF(0)=\partial_{\bar
u}F(0)=0$. 
If there are $\alpha$ and $\beta$ 
positive  with 
$T/\alpha \beta>1/4$ such that
$$
e^{|x|^2/\beta^2}\left(u_1(0)-u_2(0)\right)\, ,\,\
e^{|x|^2/\alpha^2}\left(u_1(T)-u_2(T)\right) \in L^2(\R^n), 
$$
then $u_1\equiv u_2$.

\end{theorem}

 In Theorem \ref{Theorem NL2} we did not attempt to optimize the regularity
assumption on the solutions $\,u_1,\,u_2$.

 By fixing $u_2\equiv 0$ Theorem \ref{Theorem NL2} provides a restriction on the
possible decay at two different times 
of a non-trivial solution $u_1$ of equation \eqref{e1} with
$F$ as in \eqref{F1b}. It is an open question to determine  the optimality of this
kind of result. More precisely,
for the standard semi-linear Schr\"odinger equations 
\begin{equation}
\label{NLS}
\partial_t u = i (\Delta u + |u|^{\gamma-1}u),\,\,\,\,\gamma>1,
\end{equation}
one has the  \it{standing wave  }\rm  solutions 
$$
u(x,t)=e^{\omega \,t} \varphi(x),\,\,\,\omega>0, 
$$
where $\varphi$ is the unique (up to translation) positive solution of the elliptic
problem
$$
-\Delta \varphi+\omega \varphi= |\varphi|^{\gamma-1}\varphi,
$$
which has a linear exponential decay, i.e.
$$
\varphi(x)=O(e^{-c|x|}),\,\,\,\text{as}\,\,\,|x|\to\infty,
$$
for an appropriate value of $c>0$ (see \cite{Str}, \cite{BLi}, \cite{BGK}, and
\cite{Kw}). Whether or not 
these standing waves are the solutions of \eqref{NLS} having the strongest possible decay at
two 
different times is an open question.

\vskip.04in
 Hardy's uncertainty principle also admits a formulation in terms of the heat equation
$$
\partial_tu=\Delta u,\;\;\;\;t>0,\;\;x\in\mathbb R^n,
$$
whose solution with data $\,u(x,0)=u_0(x)$ can be written as
$$
u(x,t)= e^{t \Delta}u_0(x)= \int_{\R^n} \frac{e^{-|x-y|^2/4t}}{(4\pi  t)^{n/2}}\,
u_0(y)\,dy.
$$

More precisely,  Hardy's uncertainty principle can restated in the following  equivalent forms :
\vskip.07in
\emph{ (i) If $\,u_0\in L^2(\R^n)$ and there exists  $\,T>0$ such that  $\,e^{|x|^2/(\delta^2T)}\,e^{T\Delta}u_0\in L^2(\R^n)\,$ for some $\,\delta \leq 2$, then $\,u_0\equiv 0$.}
\vskip.04in
\emph{ (ii) If $\,u_0\in \mathcal S(\R^n)$ (tempered distribution)  and  there exists  $\,T>0$ such that  $\,e^{|x|^2/(\delta^2T)}\,e^{T\Delta}u_0\in L^{\infty}(\R^n)$  for some $\,\delta< 2$, then $\,u_0\equiv 0$.
Moreover, if $\,\delta=2$, then $\,u_0$ is a constant multiple of the Dirac delta measure.}
\vskip.03in
In fact, applying Hardy's uncertainty principle to $e^{T \triangle} u_0$ one has that  $e^{\frac{|x|^2}{\delta^2 T}}e^{T \triangle} u_0$ and 
$e^{ T |\xi|^2}\widehat{e^{T \triangle} u_0}=\widehat u_0$ in $L^2(\R^n)$ with $\,2\delta\le 4$ implies $e^{\triangle}u_0\equiv 0$. Then,  backward uniqueness arguments
(see for example \cite[Chapter 3, Theorem 11]{lm60}) shows that $u_0\equiv 0$. 

In  \cite{EKPV08b} we  proved the following weaker extension of this result for parabolic operators with 
lower order variable coefficientes :
\vskip.03in	
\begin{theorem}\label{T: toremaparabolicco}
Let $u\in C([0,1] : L^2(\Rn))\cap L^2([0,T]: H^1(\Rn))$ be a solution of the IVP 
\begin{equation*}
\begin{cases}
\partial_tu=\triangle u+V(x,t)u,\ \text{in}\ \Rn\times (0,1],\\
u(x,0)=u_0(x),
\end{cases}
\end{equation*}
where 
$$
V\in L^{\infty}(\Rn\times [0,1]).
$$
If 
 $$
 u_0\;\;\;\;\text{and}\;\;\;\;e^{\frac{|x|^2}{\delta^2}}u(1)\in L^2(\Rn),
 $$
  for some $\,\delta <1$, then $u_0\equiv 0$.
\end{theorem}
\vskip.03in
		It is natural to expect  that Hardy's uncertainty principle holds in this context with 
bounded potentials $V$ and with the parameter  $\delta$ verifing the condition of the free case, i.e. $\,\delta\leq 2$.
		
Earlier results in this directions, addressing a question of Landis and Oleinik \cite{LaOl}, were obtained in \cite{EKPV06a} and  \cite{Ng}.

\vskip.03in

\section{Uncertainty Principle of Morgan type}\label{morgan}

 In \cite{Mo} G. W. Morgan proved the following uncertainty principle:

\vskip.07in
\emph{If $f(x)=O(e^{-\frac{a^p |x|^p}{p}}), \,1<p\leq 2$  and $\widehat
f(\xi)=O(e^{-\frac{(b+\epsilon)^q |\xi|^q}{q}}),\;1/p+1/q=1,\,\epsilon>0,$ 
with 
$$
ab>\Big|\cos\left(\frac{p\,\pi}{2}\right)\Big|,
$$ 
then $f\equiv 0$.}

\vskip.07in

In \cite{Ho} Beurling-H\"ormander showed :

\vskip.07in
\emph{ If
$f\in L^1(\mathbb R)$ and
\begin{equation}
\label{beurling}
\int_{\mathbb R} \int_{\mathbb R} |f(x)| |\widehat f(\xi)| e^{
|x\,\xi|}\,dx\,d\xi<\infty, \;\;\;\text{then}
\;\;\;f\equiv 0.
\end{equation}}

This result  was extended to higher dimensions $n\geq 2$ in \cite{bonamie2} and
\cite{ray} :
\vskip.07in
 \emph{If $f\in L^2(\mathbb R^n), n\geq 2$ and
\begin{equation}
\label{beurlingn}
\int_{\mathbb R^n} \int_{\mathbb R^n} |f(x)| |\widehat f(\xi)| e^{
|x\,\cdot \xi|}\,dx\,d\xi<\infty, \;\;\;\text{then}
\;\;\;f\equiv 0.
\end{equation}}
\vskip.05in

We observe that from \eqref{beurling} and \eqref{beurlingn}  it follows that :
\vskip.05in
 \emph{If 
$p\in(1,2]$, $\,1/p+1/q=1$, $\,a, \,b>0$, and 
\begin{equation}
\label{primera}
 \int_{\mathbb R^n}|f(x)|\, e^{\frac{a^p|x|^p}{p}}dx \,+
\,\int_{\mathbb R^n} |\widehat
f(\xi)| \,e^{\frac{b^q|\xi|^q}{q}}d\xi<\infty,\;\;a b\geq 1\;\Rightarrow \;f\equiv 0.
 \end{equation}}

Notice that in the case $p=q=2$ this gives us an $L^1$-version of  Hardy's uncertainty result
discussed above, and
for $p<2$ an $n$-dimensional  $L^1$-version of Morgan's uncertainty principle.

In the one-dimensional case ($n=1$), the optimal  $L^1$-version of Morgan's result
in \eqref{primera}, 
\begin{equation}
\label{primera1}
 \int_{\mathbb R}|f(x)|\, e^{\frac{a^p|x|^p}{p}}dx +
\int_{\mathbb R} |\widehat
f(\xi)| \,e^{\frac{b^q|\xi|^q}{q}}d\xi<\infty,\;\;a b>\Big|\cos\left(\frac{p\,\pi}{2}\right)\Big|\;\Rightarrow \;f\equiv 0.
 \end{equation}
was established in  \cite{bonamie2} and \cite{monki} (for further results
see \cite{bonamie1} and references therein). A sharp condition for $a,\,b,\,p$ in
\eqref{primera1} in higher dimension 
seems to be unknown.
However, in \cite{bonamie2} it was shown : 

\vskip.07in
\emph{ If $f\in
L^2(\mathbb R^n)$, $1<p\leq 2\;$ and $\,1/p+1/q=1\,$  are such that for some $j=1,..,n$,
\begin{equation}
\label{bonami11}
 \int_{\mathbb R^n}
|f(x)|e^{\frac{a^p|x_j|^p}{p}}dx<\infty\;\;+\;\;\int_{\mathbb R^n}
|\widehat f(\xi)|e^{\frac{b^q|\xi_j|^q}{q}}d\xi<\infty.
\end{equation}}

\vskip.05in
\emph{If  $a b>\Big|\cos\left(\frac{p\,\pi}{2}\right)\Big|$, then
$\;f\equiv 0$.}

\vskip.05in
\emph{If $a b<\Big|\cos\left(\frac{p\,\pi}{2}\right)\Big|$, then  there exist 
non-trivial functions satisfying \eqref{bonami11}}. 
\vskip.05in
 
Using \eqref{formula1}  the above result  can be stated in terms of 
the solution of the free Schr\"odinger
equation. In particular, \eqref{primera} can be re-written as  :

\vskip.05in
\emph{If $u_0\in L^1(\mathbb R)$ or $u_0\in L^2(\mathbb R^n)$, if $n\geq 2$,  and for
some $\,t\neq 0$
 \begin{equation}
 \label{pq}
 \int_{\mathbb R^n}\;|u_0(x)|\, e^{\frac{a^p|x|^p}{p}}dx \,+
\,\int_{\mathbb R^n}\,
|\,e^{it\Delta}u_0(x)|
\,e^{\frac{b^q|x|^q}{q(2t)^q}}dx<\infty,
\end{equation}
with 
$$
ab>\Big|\cos\left(\frac{p\,\pi}{2}\right)\Big|\;\;\;\;\text{if}\;\;\;n=1,\;\;\;\;\text{and}\;\;\;\;ab>1\;\;\;\;\text{if}\;\;\;\;\;n\geq 2,
$$
then $u_0\equiv 0$. }
\vskip.07in

Related with Morgan's uncertainty principle one has the following result due to Gel'fand and Shilov. 
In \cite{GeShi}  they considered the class $Z^p_p,\,p> 1$,
defined as  the space of all functions $\varphi(z_1,..,z_n)$ which are
analytic for all values of $z_1,..,z_n\in \mathbb C$ and such that
$$
|\varphi(z_1,..,z_n)|\leq C_0\, e^{\sum_{j=1}^n\,\epsilon_j\,C_j\,|z_j|^p},
$$
where the $C_j,\,j=0,1,..,n$ are positive constants and $\epsilon_j=1$ for
$z_j$ non-real and $\epsilon_j=-1$ for $z_j$ real, $j=1,..,n$, and showed
that the  Fourier transform of the function space $Z_p^p$ is the space
$Z_q^q$,
with $\,1/p+1/q=1$.
\vskip.05in
Notice that  the class $Z_p^p$ with $p\geq 2$ is closed with respect to
multiplication by $\,e^{i c |x|^2}$. Thus, if $u_0\in Z^p_p,\,p\geq 2$, then by
\eqref{formula1} 
one has that 
$$
|e^{it\Delta}u_0(x)|\leq d(t)\,e^{-a(t)|x|^q},
$$
 for some functions
$\,d,\,a\,:\,\mathbb R\to(0,\infty)$.
\vskip.03in
In \cite{EKPV08m}  the following results were established:

\begin{theorem}\label{Theorem 22}
Given $\,p\in(1,2)$ there exists $\,M_p>0$ such that  for any solution  $u
\in C([0,1] :L^2(\Rn))$ of
\begin{equation*}
\label{E: 1.111}\partial_tu=i\left(\triangle
u+V(x,t)u\right),\;\;\;\text{in}\;\;\;\;\;\Rn\times [0,1],
\end{equation*}
with  $V=V(x,t)$  complex valued,  bounded (i.e.
$\|V\|_{L^{\infty}(\R^n\times[0,1])}\leq C$)
and
\begin{equation}
\label{14}
\lim_{R\rightarrow +\infty}\|V\|_{L^1([0,1] : L^\infty(\Rn\setminus B_R))}=0,
\end{equation}
satisfying that for some constants $\,a_0,\,a_1,\,a_2>0$
\begin{equation}
\label{12}
\int_{\mathbb R^n} |u(x,0)|^2\,e^{2a_0 |x|^p}dx < \infty,
\end{equation}
and for any $k\in\Z^+$
\begin{equation}
\label{13}
\int_{\mathbb R^n} |u(x,1)|^2\,e^{2k |x|^p}dx < a_2 e^{2 a_1 k^{q/(q-p)}},
\end{equation}
$1/p+1/q=1$, if
\begin{equation}
\label{conditionp}
\,a_0\,a_1^{(p-2)} > M_p,
\end{equation}
then $\,u\equiv 0$.
\end{theorem}

\begin{corollary}\label{Corollary 22}
Given $\,p\in(1,2)$ there exists $N_p>0$ such that if
\newline $u\in C([0,1]:L^2(\mathbb R^n))$ is a solution of
$$
\partial_t u=i (\Delta u +V(x,t)u),
$$
with  $V=V(x,t)$  complex valued,  bounded (i.e.
$\|V\|_{L^{\infty}(\R^n\times[0,1])}\leq C$) and 
$$
\lim_{R\to\infty} \,\int_0^1\,\sup_{|x|>R} |V(x,t)| dt=0,
$$
and there exist $\,\alpha,\,\beta>0$ such that
\begin{equation}
\label{con1}
\int_{\mathbb R^n}
|u(x,0)|^2e^{2\,\alpha^p\,|x|^p/p}dx\;\,\,+\,\,\;\int_{\mathbb
R^n}
|u(x,1)|^2e^{2\,\beta^q\,|x|^q/q}dx<\infty,
\end{equation}
$\,1/p+1/q=1$, with
\begin{equation}
\label{conditionp2}
\;\alpha\,\beta > N_p,
\end{equation}
then $\;u\equiv 0$.
\end{corollary}

As a  consequence of  Corollary \ref{Corollary 22}  one obtains the
following   result concerning the uniqueness of solutions for the semi-linear
equations \eqref{e1} with $F$ as in \eqref{F1b}
\begin{equation}
\label{E: NLS}
i \partial_t u + \triangle u = F(u,\overline u).
\end{equation}

\begin{theorem}
\label{Theorem 23}

Given $\,p\in(1,2)$ there exists $\,N_p>0$ such that  if
$$
u_1,\,u_2 \in C([0,1] : H^k(\R^n)),
$$
are strong solutions of \eqref{E: NLS}  with $k\in \Z^+$, $k>n/2$,
$F:\C^2\to \C$, $F\in C^{k}$  and $F(0)=\partial_uF(0)=\partial_{\bar
u}F(0)=0$, and there exist $\,\alpha,\,\beta>0$
such that
\begin{equation}
\label{con2}
e^{\alpha^p\,|x|^p/p}\left(u_1(0)-u_2(0)\right),\;\;\;\
e^{\beta^q\,|x|^q/q}\left(u_1(1)-u_2(1)\right) \in L^2(\R^n),
\end{equation}
$1/p+1/q=1$, with
\begin{equation}
\label{conditionp2b}
\,\alpha\,\beta > N_p,
\end{equation}
then $u_1\equiv u_2$.

\end{theorem}

Notice that the conditions \eqref{conditionp} and  \eqref{conditionp2} are
independent of the size of
the potential and there is not any \it{a priori }\rm   regularity assumption  on the
potential $V(x,t)$.

The result in \cite{bonamie2}, see \eqref{bonami11}, can be extended to our setting 
with an non-optimal constant. More precisely,

\begin{corollary}\label{Corollary 24}
The conclusions in Corollary \ref{Corollary 22} still hold with a different
constant $N_p>0$ if one replaces the hypothesis  \eqref{con1} by the following one
dimensional version
\begin{equation}
\label{conn=1}
\int_{\mathbb R^n}
|u(x,0)|^2e^{2\,\alpha^p\,|x_j|^p/p}dx<\infty\,\;\,\,+\,\,\;\int_{\mathbb
R^n}
|u(x,1)|^2e^{2\,\beta^q\,|x_j|^q/q}dx<\infty,
\end{equation}
for some $j=1,..,n$. 
\end{corollary}

 Similarly, the non-linear version of Theorem \ref{Theorem 23} 
still holds, with different constant $N_p>0$,  if one replaces the hypothesis
\eqref{con2} by 
$$
e^{\alpha^p\,|x_j|^p/p}\left(u_1(0)-u_2(0)\right),\;\;\;\
e^{\beta^q\,|x_j|^q/q}\left(u_1(1)-u_2(1)\right) \in L^2(\R^n),
$$
for $j=1,..,n$. 

 In \cite{EKPV08m} we did  not attempt to give an estimate
of the universal constant $N_p$. 

The limiting case $\,p=1$ will be considered in the next section.

The main idea in the proof of these results is to combine an
upper estimate with a lower one to obtain the desired result. The upper estimate is based on the decay hypothesis on the solution at two different
times 
(see Lemma \ref{ultimo}).
In previous works we had been able to establish these estimates from  assumptions
that at time $t=0$ 
and $t=1$ involving 
the
same weight. However, in our case (Corollary \ref{Corollary 22}) we have
different weights at time $t=0$ and $t=1$. To
overcome this difficulty, we carry out the details with the weight
$e^{a_j|x|^p},\,1<p<2$, $j=0$ at $t=0$ and $j=1$ at $t=1$, with $a_0$ fixed and
$a_1=k\in\mathbb Z^+$
as in \eqref{13}. Although the powers $\,|x|^p\,$ in the exponential are
equal at time $t=0$ and $t=1$ to apply our estimate (Lemma \ref{ultimo})
we also need  to have the same constant in front of them. To achieve this we
apply the conformal or Appell transformation described above, to get  solutions and
potentials, 
whose bounds  depend on $k\in\mathbb Z^+$.
Thus we have to  consider a  family of solutions and obtain estimates on their
asymptotic value as
$k\uparrow \infty$.

The proof of the lower estimate is based on the
positivity of the commutator operator obtained by conjugating the equation
with the appropriate exponential weight, (see Lemma \ref{L: freq1} in the appendix)

\section{Paley-Wiener Theorem and Uncertainty Principle of Ingham type}\label{ingham}

 This section is concerned with the limiting case $p=1$ in the previous section.

It is easy to see that if $f\in L^1(\R^n)$ is non-zero and has compact support, then
$\,\widehat f $
cannot satisfy a condition of the type
$\widehat f(y)=O(e^{-\epsilon |y|})$ for any $\epsilon>0$.
However, it may be possible to have  $ f\in L^1(\R^n)$ a non-zero function with
compact support, such that
$\widehat f(\xi)=O(e^{-\epsilon(y) |y|})$, $\epsilon(y)$ being a positive function
tending to zero as 
$|y|\to \infty$.

In the one-dimensional case ($n=1$) soon after Hardy's result described above, A. E.
Ingham \cite{In} 
proved the following :
\vskip.05in
\emph{There exists $f\in L^1(\R)$ non-zero, even, vanishing outside an interval such
that
$\widehat f(y)=O(e^{-\epsilon(y) |y|})$ with $\epsilon(y)$ being a positive
function tending to zero at infinity
if and only if
$$
\int^{\infty} \frac{\epsilon(y)}{y}\,dy<\infty.
$$}

 In a similar direction  the Paley-Wiener Theorem \cite{PW} gives a characterization of
a function 
or distribution with compact support in term of analyticity properties of its
Fourier transform.

 Regarding our results discussed above it would be interesting  to identify a class
of potentials $V(x,t)$ for which 
a  result of the following kind holds:

\vskip.05in
If $u\in C([0,1]:L^2(\R^n))$ is a non-trivial solution of the IVP
\begin{equation}
\label{0007}
\begin{cases}
\begin{aligned}
&\partial_tu=i(\triangle u+V(x,t)u), \;\;\,\,(x,t)\in\R^n\times [0,1],\\
&u(x,0)=u_0(x),
\end{aligned}
\end{cases}
\end{equation}
with $u_0\in L^2(\R^n)$ having compact support, then $ e^{\epsilon
|x|}\,u(\cdot,t)\notin L^2(\R^n)$ 
for any $\epsilon>0$ and any $t\in(0,1]$.
\vskip.03in

In this direction we have the following result which will appear in \cite{EKPV12}:

\begin{theorem}\label{2012} Assume that 
 $u\in C([0,1]:L^2(\mathbb R^n))$ is a strong solution of the IVP \eqref{007}
with
\begin{equation}
\label{hyp1-2012}
supp\,u_0\subset B_R(0)=\{x\in\mathbb R^n\,:\,|x|\leq R\},
\end{equation}
\begin{equation}
\label{hyp2-2012}
\int_{\mathbb R^n}\,e^{2a_1|x|}\,|u(x,1)|^2\,dx<\infty,\;\;\;\;\;\;a_1>0, 
\end{equation}
and
\begin{equation}
\label{hyp3-2012}
 \|V\|_{L^{\infty}(\mathbb R^n\times [0,1])}=M_0,
\end{equation}
with 
\begin{equation}
\label{hyp4-2012}
\lim_{R\rightarrow +\infty}\|V\|_{L^1([0,1] : L^\infty(\Rn\setminus B_R))}=0.
\end{equation}
Then, there exists $b=b(n)>0$ (depending only on the dimension $n$) 
such that if
$$
\frac{a_1}{R\,(1+M_0)}\geq b,
$$
then $\,u\equiv 0$.
\end{theorem}

\vskip.05in

A similar question can be raised  for results of the type described above due to A.
E. Ingham in \cite{In}
and  possible extensions to higher dimensions $n\geq 2$.
\vskip.05in
It would be interesting to obtain extensions of the above results 
characterizing the
decay of the solution $u(x,t)$ to the equation \eqref{e1} with $F$ as in \eqref{F1b}
associated 
to data $u_0\in L^2(\R^n)$ with compact support or with $u_0\in C_0^{\infty}(\R^n)$.
In this direction, some results can be deduced as  a consequence of Theorem \ref{2012},
see \cite{EKPV12}.

\section{Hardy's Uncertainty Principle in a half-space}\label{half}

In the introduction we have briefly reviewed some uniqueness results established for
solutions 
of the Schr\"odinger equation vanishing at two different times in a semi-space of
$\,\R^n$, 
(see \cite{BZ}, \cite{ds}, \cite{IK04}, \cite{IK06}, \cite{EKPV08b}). In section 2,
we have studied
uniqueness results gotten under the hypothesis that the solution 
of the Schr\"odinger equation at two different times has an appropriate Gaussian
decay, in the whole space $\R^n$.
In this section, we shall deduce a unified result, i.e. a uniqueness result 
 under the hypothesis that at two different times the solution of the Schr\"odinger
equation has Gaussian decay in 
just a semi-space of $\,\R^n$.

\begin{theorem}\label{hardyhalf} Assume that 
 $u\in C([0,1]:L^2((0,\infty)\times \mathbb R^{n-1}))$ is a strong solution of the IVP
\begin{equation}
\begin{cases}
\begin{aligned}
\label{eq441}
&\partial_tu=i(\Delta + V(x,t))u,\\
&u(x,0)=u_0(x),
\end{aligned}
\end{cases}
\end{equation}
with
\begin{equation}
\label{extrahyp}
\int_0^1\,\int_{1/2}^{3/2}\,|\partial_{x_1}u(x,t)|^2\,dx\,dt<\infty, 
\end{equation}
\begin{equation}
\label{hyp442}
 V:\mathbb R^n\times [0,1]\to\mathbb C,\,\,\,\,\,\,\,V\in L^{\infty}(\mathbb
R^n\times [0,1]),
\end{equation}
and \begin{equation}
 \label{condition44}
\lim_{R\rightarrow +\infty}\int_0^1\|V(t)\|_{L^\infty(\{x_1>R\})}\,dt =0.
\end{equation}
Assume that 
\begin{equation}
\begin{aligned}
\label{443}
&\int_{x_1>0} \,e^{c_0\,|x_1|^2}\,|u(x,0)|^2\, dx <\infty,\\
\\
&\int_{x_1>0} \,e^{c_1\,|x_1|^2}\,|u(x,1)|^2\, dx <\infty,
\end{aligned}
\end{equation}
with $c_0,\,c_1>\,0$ sufficiently large.
Then $\,u\equiv 0$.
\end{theorem}

\underline{Remarks} : (a)   Note that in Theorem \ref{hardyhalf}, the solution does not need to be defined for $\,x_1\leq 0$.
In this sense, this is a stronger result that the uniqueness results in \cite{BZ}, \cite{KPV02}, \cite{IK04}, \cite{IK06},
and \cite{ds}, which required that the solution be defined in $\,\mathbb R^n\times [0,1]$ and be $C([0,1]:L^2(\mathbb R^n))$.

On the other hand, we need to assume the condition \eqref{extrahyp}. Note that \cite{KPV02} also needs an extra assumption on $\,\nabla u$, 
stronger that \eqref{extrahyp}, but that in \cite{IK04}, which  among other things removed any extra assumption on $\,\nabla u$, but still required
the solution to be defined in $\,\mathbb R^n\times [0,1]$ and be in $C([0,1]:L^2(\mathbb R^n))$. If in the setting of Theorem \ref{hardyhalf}
we know that $\,u\,$ is a solution in $\,\mathbb R^n\times [0,1]$ and is in $C([0,1]:L^2(\mathbb R^n))$, then we can dispose the hypothesis
\eqref{extrahyp} as follows: 

 First as in the first step of the proof of Theorem \ref{hardyhalf}, we can use the Appell transformation to reduce to the case
$c_1=c_2=2\gamma$. Then, using $\,\varphi(x_1)$ a \lq\lq regularized" convex function which agrees with $\,x_1^+$ for $x_1>1$ , $x_1<-1$, 
an application of Lemma \ref{L: freq1} and Corollary \ref{last} in the appendix yields the estimate
$$
\sup_{0\leq t\leq 1}\,\int\,e^{2\gamma(x_1^+)^2}|u(x,t)|^2dx+\int_0^1\int_{x_1>2}\,t(1-t)|\nabla u(x,t)|^2e^{2\gamma(x_1^+)^2}dxdt<\infty.
$$
Once this is obtained, by restricting our attention to 
$$
(2,\infty)\times \mathbb R^{n-1}\times [\delta,1-\delta],
$$
for each $\,\delta>0$, we are in the situation of  Theorem \ref{hardyhalf}, and hence $\,u\equiv 0$ on $\{x_1>2\}\times[0,1]$. Finally, Izakov's result
in \cite{Iza} concludes that $\,u\equiv 0$ (more precisely, the version of Izakov's result proved in \cite{IK04}, which does not require $\,\nabla u$ to exist
for $-1<x_1<1).$
\vskip.05in
(b) We have seen that Theorem \ref{hardyhalf} includes many of the uniqueness results for solutions vanishing at two different times
in a semi-space. In comparison with the results in section 2, since the extra assumption \eqref{extrahyp} can be recovered as in remark (a) when
the solution is defined in $\,\mathbb R^n\times [0,1]$ and is in $C([0,1]:L^2(\mathbb R^n))$, the only weakness is that the provide an optimal
estimate for the constants $\,c_1,\,c_2$, but on the other hand deals with solutions only defined in $(0,\infty)\times \mathbb R^{n-1}\times [0,1]$. 
\vskip.05in
(c) In Theorem  \ref{hardyhalf} the direction $\,\vec e_1$ can be replaced by any
other  $\,\omega\in\mathcal S^{n-1}$.
\vskip.1in

\underline{Proof of Theorem \ref{hardyhalf}}:  The strategy of the proof follows closely the one in \cite{EKPV06}. We divide the proof into three steps.

\vskip.05in

\underline{First Step } : Reduction to the case $c_0=c_1=2 \gamma$.

\vskip.05in

 This follows by using the conformal or Appell transformation introduced in section 2
(see \eqref{2.1}-\eqref{externalforce}), combined with the observation that the set $\{x_1>0\}$ remains invariant. 

\vskip.05in

\underline{Second Step } : Upper Bounds.

\vskip.05in

We define
$$
v(x,t)=\theta(x_1)\,u(x,t),
$$
with $\,\theta \in C^{\infty}(\R)$, non-decreasing with $\,\theta(x_1)\equiv 1\,$ if
$\,x_1>3/2$, and $\theta(x_1)\equiv 0\,$ if $\,x_1<1/2$.
Therefore,
\begin{equation}
 \label{FFF}
\partial_t v=i\,\Delta v + i\,V(x,t) v +
i\,F(x,t),\;\;\;\,\,\,\;\;\;F(x,t)=2\,\partial_{x_1}u\,\theta'(x_1)+
u\,\theta''(x_1).
\end{equation}

Using \eqref{extrahyp} we can apply Lemma \ref{ultimo} to get that
\begin{equation}
\begin{aligned}
\label{uno44}
&\sup_{0\leq t\leq 1}\| e^{\lambda\cdot x_1} v(\,\cdot\,,t)\|_{L^2(\mathbb \R^n)} \\
&\leq c_n
\Big(\|e^{\lambda\cdot x_1} v(0)\|_{L^2(\mathbb \R^n)} + \|e^{\lambda\cdot x_1}
v(1)\|_{L^2(\mathbb \R^n)}\\
& +\int_0^1
\|e^{\lambda\cdot x_1}\, F(\cdot, t)\|_{L^2(\mathbb \R^n)} dt
+ \int_0^1
\|e^{\lambda\cdot x_1}\, V\,\chi_{\{x_1<R\}}v(\cdot, t)\|_{L^2(\mathbb \R^n)} dt\Big),
\end{aligned}
\end{equation}
for some fixed $R$ sufficiently large. Thus, using \eqref{extrahyp} 
\begin{equation}
\begin{aligned}
\label{uno444}
&\sup_{0\leq t\leq 1}\| e^{\lambda\cdot x_1} v(\,\cdot\,,t)\|_{L^2(\mathbb \R^n)} \\
&\leq c_n
\Big(\|e^{\lambda\cdot x_1} v(0)\|_{L^2(\mathbb \R^n)} + \|e^{\lambda\cdot x_1}
v(1)\|_{L^2(\mathbb \R^n)}\\
& + c\,e^{c\,|\lambda|}
+ c\,\|V\|_{\infty} \,e^{c\,|\lambda|\,R}\Big).
\end{aligned}
\end{equation}
Thus, from the formula \eqref{est1} (with $p=2$ and $n=1$) and \eqref{uno444} 
we obtain that
$$
\aligned
&\sup_{0\leq t\leq 1}\| e^{\gamma\,|x_1|^2} v(\,\cdot\,,t)\|_{L^2(\mathbb \R^n)}
\\
&
\,\,\,\leq 
\Big(\| e^{\gamma\,|x_1|^2} v(0)\|_{L^2(\mathbb \R^n)} + \| e^{\gamma\,|x_1|^2}
v(1)\|_{L^2(\mathbb \R^n)}
+c +\,\|V\|_{\infty} \,e^{c\,\gamma\,R^2}\Big).
\endaligned
$$
Thus,
\begin{equation}
\label{step2a}
 \sup_{0\leq t\leq 1}\| e^{\gamma\,|x_1|^2} v(\,\cdot\,,t)\|_{L^2(\mathbb \R^n)}\leq
c_{\gamma}.
\end{equation}
 Combining this and the equation for $\,v\,$ we shall get a smoothing estimate.
Using the notation 
$$
H(t)=\|f\|^2_{L^2(\R^n)}=\|f\|^2,
$$
with 
$$
f(x,t)= e^{\gamma|x_1|^2}\,v(x,t)
$$
and the abstract Lemma \ref{L: freq1} (see the appendix) one formally has that
\begin{equation}
\label{upper-smooth}
\begin{aligned}
\partial_t^2H &\leq  2\partial_t\text{\it Re}\left(\partial_tf-\mathcal Sf-\mathcal
Af,f\right)\\
&+ 2\left(\mathcal S_tf+\left[\mathcal S,\mathcal A\right]f,f\right) + \|\,e^{\gamma
|x_1|^2}(F + V\,v)\|^2,
\end{aligned}
\end{equation}
with
$$
e^{\gamma|x_1|^2}(\partial_t-i\,\Delta) (e^{-\gamma|x_1|^2}f) = \partial_t f 
-\mathcal Sf-\mathcal Af= e^{\gamma |x_1|^2}(F + V v),
$$
where  $\, \mathcal S = - i \gamma (4x_1\,\partial_{x_1}+2)$ is symmetric,
$\mathcal A=i(\Delta + 4\gamma x_1^2)$ is skew-symmetric, and $\,F\,$ as in
\eqref{FFF}.
Since,
$$
[\mathcal S ; \mathcal A] = -8 \gamma \partial_{x_1}^2+ 16 \gamma^2\,x_1^2.
$$
using the inequality
$$
\aligned
&\int_{\R^n}\,(|\partial_{x_1}f|^2+4\gamma^2|x_1|^2|f|^2)\,dx =
\int_{\R^n}\,e^{2\,\gamma|x_1|^2}\,(|\partial_{x_1}u|^2-2\gamma\,|u|^2)dx\\
&\geq 2\,\gamma\,\int_{\R^n}\,|f|^2\,dx.
\endaligned
$$
together with Corollary \ref{last} we conclude that
\begin{equation}
 \label{009}
\int_0^1\,\int \,t(1-t) \,|\partial_{x_1}
v(x,t)|^2\,e^{2\,\gamma|x_1|^2}\,e^{2\,\gamma |x_1|^2}\,dx\,dt \leq c_{\gamma}.
\end{equation}

Combining and \eqref{step2a} and  \eqref{009} one gets that
\begin{equation}
\label{00step2}
\begin{aligned}
 &\sup_{0\leq t\leq 1}\| e^{\gamma\,|x_1|^2} v(\,\cdot\,,t)\|_{L^2(\mathbb \R^n)}\\
&
+ \int_0^1 \int t(1-t) |\partial_{x_1} v(x,t)|^2\,e^{2\,\gamma|x_1|^2} e^{2\,\gamma
|x_1|^2}| dx dt\leq c_{\gamma}.
\end{aligned}\end{equation}

\underline{Step3}

We recall the following result which is a slight variation of that proven in detail
in \cite{EKPV06} (Lemma 3.1, page 1818)  :

\begin{lemma}\label{CPDE} 
Assume that $ R>0$ and $\,\varphi : [0,1] \to \R$ is a smooth function. Then, there
exists 
$\,c=c(n;\|\varphi'\|_{\infty}+\|\varphi''\|_{\infty})>0$ such that the inequality
\begin{equation}
\label{cpde1} 
\frac{\alpha^{3/2}}{R^2}\,\Big\|\,e^{\alpha |\frac{x_1-x_{0_1}}{R}+\varphi(t)|^2}g
\Big\|_{L^2(dxdt)}
\leq c\, \Big\|\,e^{\alpha |\frac{x_1-x_{0_1}}{R}+\varphi(t)|^2}(i
\partial_t+\Delta) g \Big\|_{L^2(dxdt)}
\end{equation}
holds when $\,\alpha > c R^2 \,$ and $\,g\in C^{\infty}_0(\R^{n+1})\,$ is supported in
the set
$$
\{(x,t)=(x_1,..,x_n,t)\,\in\R^{n+1}\,:\, |\frac{x_1-x_{0_1}}{R}+\varphi(t)|\geq 1\}.
$$
\end{lemma}
\vspace{0,1 cm}
Now, we will chose $\,x_{0_1}=R/2$, $\;0\leq \varphi(t)\leq a,$ with $\,a=3/2-1/R$,
$\,\varphi(t)=a,$ on $\,3/8\leq t\leq 5/8$,
$\,\varphi(t)=0, $ for $\,t\in [0,1/4]\cup [3/4,1]$, and $\,\theta_R\in
C^{\infty}(\R)$ with $\,\theta_R(x_1)=1\,$ on $\,1<x_1<R-1$, 
and $\,\theta_R(x_1)=0\,$ for $\,x_1<1/2$ or $\,x_1>R$. 

Also we chose $\,\eta \in C^{\infty}(\R)$ with $\,\eta(x_1)=0,\;x_1\leq 1$ and 
$\,\eta(x_1)=1,\;x_1\geq 1+1/2R$.

We notice that up to translation we can assume that
\begin{equation}
\label{b}
\int_{3/8}^{5/8} \,\int_{2<x_1<3} \,|u(x,t)|^2dx dt=b\neq 0,
\end{equation}
otherwise we would have
$$
u(x,t)=0\;\;\,\;\;\,\;\text{on}\,\;\;\,\;(x,t)\;\,\,s.t.\,\,\,(x_1,t)\in
(0,\infty)\times (3/8,5/8),
$$
and thus by Izakov's result \cite{Iza} we would get that $u\equiv 0$. 

We let
\begin{equation}
 \label{defg}
g(x,t)=\theta_R(x_1)\,\eta\Big(\frac{x_1-R/2}{R}+\varphi(t)\Big)\,u(x,t).
\end{equation}
It is easy to see that $\,g$ is supported on the set
\begin{equation}
  \label{domain}
\{(x,t)\in\R^{n+1}\,:\, 1/2<x_1<R,\,\, 1/32<t<31/32,\;\;|\frac{x_1-R/2}{R}+\varphi(t)|\geq 1\}.
\end{equation}
so satisfies the hypothesis of Lemma \ref{CPDE}. Also   if $(x_1,t)\in (2,3)\times
(3/8,5/8)$ one has 
$\varphi=a$, $\,\eta\Big(\frac{x_1-R/2}{R}+a\Big)=1$ and $\theta_R=1$, hence in this domain  
$$
g(x,t)=u(x,t).
$$ 
Thus, from \eqref{domain} 
it follows that
$$
|\frac{x_1-R/2}{R}+\varphi(t)|\geq 1 + 1/R,
$$
so we have  the lower bound of \eqref{cpde1} 
$$
\frac{\alpha^{3/2}}{R^2} \,b\,e^{\alpha(1+1/R)^2},
$$
with $\,b\,$ as in \eqref{b}.
Now we shall estimate the right hand side of \eqref{cpde1}. Thus, 
\begin{equation}
\label{07}
\begin{aligned}
&(i\partial_t-\Delta)g= - \theta_R(x_1)\eta\Big(\frac{x_1-R/2}{R}+\varphi(t)\Big)
V(x,t) u(x,t)\\
&\;\;\;+\eta\Big(\frac{x_1-R/2}{R}+\varphi(t)\Big)(2\theta'(x_1)\,\partial_{x_1}u+u\,\theta_R''(x_1))\\
&\;\;\;+(i\eta'(\cdot)\,\varphi'(t)+\eta''(\cdot)\,\frac{1}{R^2}) \theta_R(x_1)
u(x,t)\equiv E_1+E_2+E_3.
\end{aligned}
\end{equation}

Choosing  $R>>\|V\|_{\infty}$, and recalling the fact that $\,\alpha>c R^2\,$ we see
that the contribution of the term 
$E_1$ involving the  potential $V$ can be absorbed by the term in the left hand side
of \eqref{cpde1}.

Next, we notice that the terms in $E_2$ involve derivatives of $\,\theta_R$
($\theta_R'$ or $\theta_R''$) so they are supported in 
the $(x,t)\in \R^n\times [0,1]$ such that
$$
1/2<x_1<1,\,\,\,\,\,\text{or}\,\,\,\,\,R-1<x_1<R.
$$
But, if $1/2<x_1<1$, it follows that
$$
\frac{x_1-R/2}{R}+\varphi(t)\leq
1/R-1/2+3/2-1/R=1,\,\,\,\,\text{so}\,\,\,\,\eta\Big(\frac{x_1-R/2}{R}+\varphi(t)\Big)=0.
$$
Thus, we only get contribution from the $(x,t)\in \R^n\times [0,1]$ such that
$\,R-1<x_1<R$, which can be bounded by 
$$
c\,\int_{1/32}^{31/32}\,\int_{R-1<x_1<R}\,(|u|^2+|\partial_{x_1}u|^2)(x,t)\,e^{\alpha(2-1/R)^2}\,dx\,dt.
$$

 Finally, we look at the contribution of the term in $E_3$ in \eqref{07}. In those
the derivatives fall on $\,\eta$, thus they
are supported in the region 
$$
1\leq \frac{x_1-R/2}{R}+\varphi(t)\leq 1 +
\frac{1}{2R},\;\;\,\,\,\,\,\frac{1}{2}<x_1<R,\,\,\,\,\,\;\;\frac{1}{32}<t<\frac{31}{32}.
$$
Hence, their contribution in \eqref{cpde1} is bounded by 
$$
c\,\int_{1/32}^{31/32}\,\int_{1/2<x_1<R}\,|u(x,t)|^2\,e^{\alpha(1+1/(2R))^2}\,dx\,dt
\leq c_{\gamma}\,e^{\alpha(1+1/(2R))^2}.
$$

Defining 
\begin{equation}
 \label{defdelta}
\delta(R)=\int_{1/32}^{31/32}\,\int_{R-1<x_1<R}\,(|u|^2+|\partial_{x_1}u|^2)(x,t)\,dx\,dt,
\end{equation}
and collecting the above information using that $\,\alpha=c_n\,R^2$ we get
$$
c\,R\,b\,e^{\alpha(1+1/R)^2}\leq c\,\delta(R)\,e^{\alpha(2-1/R)^2}+
\,c_{\gamma}\,e^{\alpha(1+1/(2R))^2}.
$$

Therefore, for $\,R\,$ sufficiently large it follows that (since $\,b\neq 0$)
$$
c\,R\,b\,e^{\alpha(1+1/R)^2}\leq c\,\delta(R)\,e^{\alpha(2-1/R)^2},
$$
and since $\,\alpha=c_n\,R^2$ one has that
$$
\delta(R)\geq b\,e^{-c_n R^2}.
$$

 To conclude we recall that the upper bounds in \eqref{00step2} gave us
$$
\delta(R)\leq c\,e^{-\gamma R^2},
$$
hence if $\gamma>c_n/2$ we conclude that $\,b=0$, which yields the desired result
$\,u\equiv 0$.

\section{Appendix}\label{aaa}

Above we have used the following abstract results established in  \cite{EKPV08b}:
\begin{lemma}\label{L: freq1}
Let $\mathcal S$ be a symmetric operator, $\mathcal A$ be a skew-symmetric one, both
 allowed to depend on the time variable. Let $G$ be
 a positive function, $f(x,t)$  a reasonable function, 
\begin{equation*}
\aligned
&\;H(t)=\left( f, f\right)=\|f\|^2_{L^2(\R^n)}=\|f\|^2\ ,\,\,\,\,\,\ D(t)=\left(
\mathcal Sf, f\right),\\
&\; \partial_t\mathcal S=\mathcal S_t
\quad \,\,\,\,\,\text{and}\,\,\,\,\,\quad N(t)=\frac{D(t)}{H(t)}\ .
\endaligned
\end{equation*}
Then,
\begin{multline}
\label{E: derivadasegunda}
\begin{aligned}
\partial_t^2H &= 2\partial_t\text{\it Re}\left(\partial_tf-\mathcal Sf-\mathcal
Af,f\right)+
2\left(\mathcal S_tf+\left[\mathcal S,\mathcal A\right]f,f\right)\\
&+\|\partial_tf-\mathcal Af+\mathcal Sf\|^2-
\|\partial_tf-\mathcal Af-\mathcal Sf\|^2
\end{aligned}
\end{multline}
and
\begin{equation*}
\dot N(t)\ge \left(\mathcal S_tf +\left[\mathcal S,\mathcal A\right]f, f\right)/H-
\|\partial_tf-\mathcal Af-\mathcal Sf\|^2/\left(2H\right).
\end{equation*}
Moreover, if
\begin{equation}\label{E: condicionesbase}
|\partial_tf-\mathcal Af-\mathcal Sf|\le M_1|f| +G,\ \text{in}\ \Rn\times
[0,1],\quad  \mathcal S_t+\left[\mathcal S,\mathcal A\right]\ge -M_0,
\end{equation}
and 
\[M_2=\sup_{[0,1]}{\|G(t)\|/\|f(t)\|}\]
is finite, then 
$\log H(t)$ is  \lq\lq logarithmically convex\rq\rq\ in $[0,1]$ and there is a
universal constant $N$ such that
\begin{equation}\label{E: convexidadlogaritmica}
H(t)\le e^{N\left(M_0+M_1+M_2+M_1^2+M_2^2\right)}H(0)^{1-t}H(1)^t,\ \text{when}\
0\le t\le 1.
\end{equation}
\end{lemma}

By multiplying the formula \eqref{E: derivadasegunda} by $\,t(1-t)$, integrating the
result over $[0,1]$  and using integration by parts, 
one gets the following \lq\lq smoothing" inequality
\begin{corollary}\label{last}
 With the same hypotheses and notation as in Lemma \ref{L: freq1}  
\begin{equation}
 \label{lastformula}
\begin{aligned}
 &2\int_0^1\,t(1-t)\left(\mathcal S_tf+\left[\mathcal S,\mathcal
A\right]f,f\right)\,dt+\int_0^1 \,H(t)\,dt\leq H(0)+H(1)\\
& + 2\int_0^1\,(1-2t)\text{\it Re}\left(\partial_tf-\mathcal Sf-\mathcal
Af,f\right)\,dt\\
&+
\int_0^1\,t(1-t)\|\partial_tf-\mathcal Af-\mathcal Sf\|^2_2\,dt.
\end{aligned}
\end{equation}
 \end{corollary}



\begin{thebibliography}{99}

\bibitem{AKS}  N. Aronszajn, A. Krzywicki, and J. Szarski, 
\emph{Unique continuation theorem for exterior differential forms on Riemannian manifolds},
Ark. Math. {\bf 4}  (1962) 417--453.

\bibitem{monki}  S. Ben Farah, and K. Mokni, 
\emph{Uncertainty principles and the $L^p-L^q$-version of Morgan's theorem on some
groups},
Russian J. Math. Physics {\bf 10}  (2003) 245--260.

\bibitem{BLi}  H. Berestycki, and P.-L. Lions, 
\emph{Nonlinear scalar field equations},
Arch. Rational Mech. Anal. {\bf 82}  (1983) 313--375.

\bibitem{BGK}  H. Berestycki, T. Gallou\"et, and O. Kavian, 
\emph{\'Equations de champs scalaires Euclidiens non lin\'eaires dans de plan},
C. R. Acad. Sci. Paris, ser I Math {\bf 297}  (1983) 307--310.

\bibitem{bonamie1}  A. Bonami, and B. Demange, \emph{A survey on uncertainty
principles related to quadratic forms}, Collect. 
Math. Vol. Extra (2006) 1--36. 

\bibitem{bonamie2}  A. Bonami, B. Demange, and P. Jaming, 
\emph{Hermite functions and uncertainty principles for the Fourier and 
the windowed Fourier transforms}, Rev. Mat. Iberoamericana {\bf 19}   (2006) 23--55.

\bibitem{BoKe} J. Bourgain, and C. E. Kenig, \emph{On  localization in the
continuous Anderson-Bernoulli model in higher dimensions,} Invent. Math.
{\bf 161}, 2 (2005) 1432--1297.

  \bibitem{Cal} A. P. Calder\'on, \emph{Uniqueness in the Cauchy problem for partial differential equations}, Amer. J. Math., {\bf 80} (1958),
16--36.

\bibitem{Car} T. Carleman, \emph{Sur un probl\'eme d'unicit\'e pour les syst\'emes d' equations aux deriv\'ees partielles \'a deux variables
ind\'ependantes}, Ark. Math.
{\bf 26B},  (1939) 1--9.


  \bibitem{Cha} S. Chanillo, \emph{Uniqueness of solutions to Schr\"odinger equations on complex semi-sin=mple Lie groups}, Proc. Indian Acad. Sci. 
  Math. Sci. {\bf 117} (2007),
325--331.

\bibitem{CEKPV} M. Cowling, L. Escauriaza, C. E. Kenig, G. Ponce, and L. Vega,  
\emph{The Hardy Uncertainty Principle Revisited}, to appear in Indiana  
U. Math. J.

\bibitem{CoPr} M. Cowling, and J. F. Price, 
\emph{Generalizations of Heisenberg's inequality}, Harmonic Analysis (Cortona, 1982) 
Lecture Notes in Math.,{\bf 992} (1983), 443-449, Springer, Berlin.

\bibitem{CoPr2} M. Cowling, and J. F. Price,
\emph{Bandwidth versus time concentration: the Heisenberg--Pauli--Weyl inequality}, 
SIAM J. Math. Anal. {\bf 15} (1984) 151--165.

\bibitem{cruz} J. Cruz-Sanpedro, 
\emph{Unique continuation at infinity of solutions to Schr\"odinger equations with
complex potentials},
Proc. Roy. Soc. Edinburgh. {\bf 42} (1999) 143--153.

\bibitem{ds} H. Dong, and W. Staubach.  
\emph{Unique continuation for the Schr\"odinger equation with gradient vector
potentials}, 
Proc. AMS
 {\bf 135}, 7 (2007) 2141--2149.

\bibitem{EKPV06} L. Escauriaza, C. E. Kenig, G. Ponce, and L. Vega,  
\emph{On Uniqueness  Properties of Solutions of Schr\"odinger  Equations,} 
Comm. PDE. {\bf 31}, 12 (2006) 1811--1823.

\bibitem{EKPV07} L. Escauriaza, C. E. Kenig, G. Ponce, and L. Vega,  
\emph{On  Uniqueness Properties of Solutions of the k-generalized KdV,} 
J. of Funct. Anal. {\bf 244}, 2 (2007) 504--535.

\bibitem{EKPV06a} L. Escauriaza, C. E. Kenig, G. Ponce, L. Vega, \emph{Decay at
Infinity of Caloric Functions within Characteristic Hyperplanes}, Math. Res. Letters
{\bf 13}, 3 (2006) 441--453.

\bibitem{EKPV08a} L. Escauriaza, C. E. Kenig, G. Ponce, and L. Vega,  
\emph{Convexity of Free Solutions of Schr\"odinger Equations with Gaussian Decay}, 
Math. Res. Lett. {\bf 15}, 5 (2008) 957--971.

\bibitem{EKPV08b} L. Escauriaza, C. E. Kenig, G. Ponce, and L. Vega,  
\emph{Hardy's Uncertainty Principle, Convexity and Schr\"odinger Evolutions}, 
J. European Math. Soc. {\bf 10}, 4 (2008) 883--907.

\bibitem{EKPV08m} L. Escauriaza, C. E. Kenig, G. Ponce, and L. Vega,  
\emph{Uncertainty principle of Morgan type and Schr\"odinger evolution}, J. London
Math. Soc.
{\bf 81}, (2011) 187--207.

\bibitem{EKPV09} L. Escauriaza, C. E. Kenig, G. Ponce, and L. Vega,
\emph{The sharp Hardy Uncertainty Principle for Schr\"odinger evolutions}, Duke Math. J. 
{\bf 155}, (2010)  163--187.

\bibitem{EKPV10} L. Escauriaza, C. E. Kenig, G. Ponce, and L. Vega, 
\emph{Unique continuatuion for Schr\"odinger evolutions, with 
applications to profiles of concentration and traveling waves},  
Comm. Math. Phys.  {\bf 305}, (2011) 487--512.

\bibitem{EKPV12} L. Escauriaza, C. E. Kenig, G. Ponce, L. Vega, \emph{A Theorem of Paley-Wiener Type for  Schr\"odinger Evolutions}, to appear.

\bibitem{EsSS} L. Escauriaza, G. Seregin, and V. \v Sver\'ak \emph{$L^{3,\infty}$ solutions to the Navier-Stokes equations and backward uniqueness}, Russ. Math. Surv. 
{\bf 58}, (2003) 211-250.

\bibitem{EsVe} L. Escauriaza,  and S. Vessella \emph{Optimal three cylinder inequalities for solutions to parabolic equations with Lipschitz leading coefficients}, Contemp. Math.  
{\bf 333}, Inverse Problems theory and applications (2002) 79--87.

\bibitem{GeShi} I. M. Gel'fand, and G. E. Shilov, \emph{Fourier transforms of
rapidly increasing functions and questions of uniqueness of the solution
of Cauchy's problem}, Uspehi Matem. Nauk {\bf 8}, (1953), 3--54.

\bibitem{Had} J. Hadamard, \emph{Le probl\'ome de Cauchy et les \'equations aux deriv\'ees partielles lin\'eaires hyperboliques}, 
Hermann, Paris  (1932) 

\bibitem{Hardy} G. H. Hardy, \emph{A Theorem Concerning Fourier Transforms}, J.
London Math. Soc. s1-8 (1933) 227--231.


\bibitem{HKT1} N. Hayashi, K. Nakamitsu, and M. Tsutsumi, 
\emph{On solutions of the initial value problem for the nonlinear Schr\"odinger
equations in one space dimension},  
Math. Z. {\bf 192} (1986) 637--650.

\bibitem{HKT2} N. Hayashi, K. Nakamitsu, and M. Tsutsumi, \emph{On solutions of the
initial
value problem for the nonlinear Schr\"odinger equations}, J. Funct. Anal. {\bf 71}
(1987) 218--245.

\bibitem{Ho} L. H\"ormander, \emph{A uniqueness theorem of Beurling for Fourier
transform pairs}, Ark. Mat. {\bf 29}, 2 (1991) 237--240.

\bibitem{In} A. E. Ingham, 
\emph{A note on Fourier transforms}  J. London Math. Soc. s1-9 (1934) 29--32.



\bibitem{IK04} A. D. Ionescu, and C. E. Kenig, 
\emph{$L^p$-Carleman inequalities and uniqueness of solutions of nonlinear
Schr\"odinger equations,} 
Acta Math. {\bf 193}, 2 (2004) 193--239.

\bibitem{IK06} A. D. Ionescu, and C. E. Kenig, 
\emph{Uniqueness properties of solutions of Schr\"odinger equations,} J. Funct. Anal. 
{\bf 232} (2006) 90--136.

\bibitem{Iza} V. Izakov
\emph{Carleman type estimates in an anisotropic case and applications,} J. Diff . Eqs. 
{\bf 105} (1993) 217--238.

\bibitem{Ke1} C. E. Kenig,
\emph{Some recent quantitative unique continuation theorem}, Rendiconti Accademia Nazionale dell Scienze, vol. {\bf XXIX} (2005) 231--242.

\bibitem{Ke2} C. E. Kenig, 
\emph{Some recent applications of unique continuations}, Contemp. Math. 
{\bf 439} (2007) 25--56.

\bibitem{KeMe} C. E. Kenig, and F. Merle
\emph{Global well-posedness scattering and blow-up for the energy-critical focusing non-linear wave equation}, Acta Math. 
{\bf 201} (2008) 147--212.

\bibitem{KPV02} C. E. Kenig, G. Ponce, and L. Vega,
\emph{On unique continuation for nonlinear Schr\"odinger equations}, Comm. Pure
Appl. Math. 
{\bf 60} (2002) 1247--1262.


\bibitem{KoLa} V. A. Kondratiev, and E. M. Landis
\emph{Qualitative properties of the solutions of a second-order nonlinear equation}, Encyclopedia of Math. Sci. 32 
(Partial Differential equations III) Springer-Verlag, Berlin 
{\bf } (1988).



\bibitem{Kw} K. M. Kwong, 
\emph{Uniqueness of positive solutions of $\Delta u-u+u=0$ in $\R^n$},
Arch. Rational Mech. Anal. {\bf 105}  (1989) 243-266.

\bibitem{LaOl}   E. M. Landis, and O. A. Oleinik, \emph{Generalized analyticity and certain properties, of solutions of elliptic and parabolic equations, that are connected with it}, 
Uspehi Mat. Nauk {\bf 29} (1974), 190-206.

\bibitem{lm60}  J. L. Lions, and  B. Malgrange, \emph{Sur l'unicit\'e r\'etrograde dans les probl\`emes mixtes paraboliques,} Math. Scan. {\bf8} (1960) 277--286.


\bibitem{Me} V. Z. Meshkov, \emph{On the possible rate of decay at infinity
of solutions of second-order partial differential equations}, Math. USSR Sbornik  
{\bf 72} (1992), 343--361.

\bibitem{Mo} G. W. Morgan, \emph{A note on Fourier transforms}, J. London Math. Soc.
  {\bf 9} (1934),   187-192.
  
  \bibitem{Mu} C. M\"uller, \emph{On the behavior of the solution of the differential equation $\Delta u=f(x,u)$ in the neighborhood of a point}, 
Comm. Math. Pure Appl. 
  {\bf 1} (1954),   505--515.
  
  
  \bibitem{ray}  E. K. Naranayan, and S. K. Ray, \emph{Beurling's theorem in $\mathbb
R^n$}, preprint.

 \bibitem{Ng} T. Nguyen, \emph{On a question of Landis and Oleinik}, Trans. Amer. Math. Soc. {\bf 362} (2010) 2875--2899

\bibitem{PW} R. Paley, and N. Wiener, \emph{Fourier transform in the complex domain},
Amer. Math. Soc  Providence RI (1934).




\bibitem{SaSc} J.-C. Saut, and B. Scheurer, 
\emph {Unique continuation for some evolution equations},
J. Diff. Eqs. {\bf  66}
(1987), 118--139.


\bibitem{Se}  G. Seregin,  \emph{A certain necessary condition of potential blow-up for Navier-Stokes equations}, preprint, arXiv: 1104.3615

\bibitem{SST} A. Sitaram, M. Sundari, and S. Thangavelu, \emph{Uncertainty
principles on certain Lie groups},  
Proc. Indian Acad. Sci. Math. Sci. 
{\bf 105} (1995), 135-151.


\bibitem{Str} W. A. Strauss, \emph{Existence of solitary waves in higher dimensions},  
Comm. Math. Phys. 
{\bf 55} (1977), 149--162.

\bibitem{Ve}  S. Vessella \emph{Three cylinder inequalities and unique continuation properties for parabolic equations},  Atti Accad. Naz. Lincei Cl. Sci. Fis. Mat. Natur. Rend. Lincei (9) Mat. Appl. {\bf 13} (2002), no. 2, 107Ð120


\bibitem{BZ} B. Y. Zhang, \emph{Unique continuation properties of the 
nonlinear Schr\"odinger equations}, Proc. Roy. Soc. Edinburgh. {\bf 127} (1997)
191--205.

\end{thebibliography}
\end{document}